\documentclass[11pt]{amsart}

\usepackage{geometry}
\usepackage{todonotes}

\usepackage{amsfonts, amssymb, amscd}
\numberwithin{equation}{section}

\usepackage[symbol]{footmisc}

\usepackage{bm}
\usepackage{verbatim}
\usepackage{mathrsfs}
\usepackage{graphicx}
\usepackage{tikz-cd}
\usepackage{subcaption}
\usepackage{listings}
\usepackage{subfiles}
\usepackage[toc,page]{appendix}
\usepackage{mathtools}
\usepackage{comment}
\usepackage{enumerate}
\usepackage{enumitem}
\usepackage[all]{xy}

\usepackage{graphicx}
\graphicspath{{images/}}

\usepackage{appendix}
\usepackage{hyperref}
\hypersetup{
    colorlinks=true,
    citecolor=red,
    linkcolor=blue,
    filecolor=magenta,      
    urlcolor=red,
}
\newcommand{\Qq}{\mathbb{Q}}

\newcommand{\Rr}{\mathbb{R}}

\newcommand{\Center}{\operatorname{center}}

\newcommand{\Exc}{\operatorname{Exc}}

\newcommand{\ct}{\operatorname{ct}}

\newcommand{\mld}{{\operatorname{mld}}}
\newcommand{\pld}{{\operatorname{pld}}}
\newcommand{\tmld}{{\operatorname{tmld}}}

\newcommand{\lct}{\operatorname{lct}}

\newcommand{\Supp}{\operatorname{Supp}}

\newcommand{\mult}{\operatorname{mult}}

\newcommand{\Ii}{\Gamma}

\newtheorem{thm}{Theorem}[section]

\newtheorem{lem}[thm]{Lemma}

\newtheorem{claim}[thm]{Claim}

\theoremstyle{definition}
\newtheorem{defn}[thm]{Definition}

\theoremstyle{definition}

\newtheorem{ex}[thm]{Example}

\theoremstyle{definition}

\begin{document}

\title{Infinitesimal structure of log canonical thresholds}
\author{Jihao Liu, Fanjun Meng, and Lingyao Xie}

\date{\today}

\begin{abstract}
We show that log canonical thresholds of fixed dimension are standardized. More precisely, we show that any sequence of log canonical thresholds in fixed dimension $d$ accumulates either i) in a way which is similar to how standard and hyperstandard sets accumulate, or ii) to log canonical thresholds in dimension $\leq d-2$. This provides an accurate description on the infinitesimal structure of the set of log canonical thresholds. We also discuss similar behaviors of minimal log discrepancies, canonical thresholds, and K-semistable thresholds.
\end{abstract}

\address{Department of Mathematics, Northwestern University, 2033 Sheridan Road, Evanston, IL 60208, USA}
\email{jliu@northwestern.edu}

\address{Department of Mathematics, Johns Hopkins University, 3400 N. Charles Street, Baltimore, MD 21218, USA}
\email{fmeng3@jhu.edu}

\address{Department of Mathematics, The University of Utah,  155 South 1400 East, JWB 233, Salt Lake City, UT 84112, USA}
\email{lingyao@math.utah.edu}

\maketitle

\tableofcontents

\section{Introduction}

We work over the field of complex numbers $\mathbb C$. For any set $\Ii\subset\mathbb R$, we let $\partial\Ii$ be the set of accumulation points of $\Ii$ and $\bar\Ii:=\Ii\cup\partial\Ii$ the closure of $\Ii$. We let $\partial^0\Ii:=\bar\Ii$, and denote the set of $k$-th order accumulation points of $\Ii$ by $\partial^k\Ii$ for any $k>0$. It is clear that $\partial^k\Ii=\partial^k\bar\Ii$ for any non-negative integer $k$.

\medskip

\noindent\textbf{Log canonical thresholds}. The log canonical threshold (lct for short) is a fundamental invariant in algebraic geometry. It originates from analysis which measures the integrability of a holomorphic function. In birational geometry, the log canonical threshold measures the complexity of the singularities of a triple $(X,B;D)$ where $(X,B)$ is a pair and $D$ is an effective $\Rr$-Cartier $\Rr$-divisor.

\begin{defn}
Let $(X,B)$ be a pair and $D\geq 0$ an $\Rr$-Cartier $\mathbb R$-divisor. We define
$$\lct(X,B;D):=\sup\{t\geq 0\mid (X,B+tD)\text{ is log canonical (lc)}\}$$
to be the lct of $D$ with respect to $(X,B)$. The set of lcts in dimension $d$ is defined as
\begin{align*}
    \lct(d):=\Bigg\{\lct(X,B;D)\Biggm|
    \begin{array}{r@{}l}
     \dim X=d, (X,B)\text{ is lc, } B \text{ is an effective Weil divisor},\\
      \text{ and }D\text{ is an effective }\Qq\text{-Cartier Weil divisor}\\ 
    \end{array}\Bigg\}.
    \end{align*}
\end{defn}

It is well-known that the set of log canonical thresholds of fixed dimension satisfies the ascending chain condition (ACC) \cite[Theorem 1.1]{HMX14} and their accumulation points are log canonical thresholds from lower dimension \cite[Theorem 1.11]{HMX14}. The purpose of this paper is to discuss how lcts approach their accumulation points. More precisely, we show that the sets of lcts of fixed dimension are \emph{standardized sets}. Roughly speaking, this says that the infinitesimal behavior of the sets of lcts is similar to the behavior of standard and hyperstandard sets, especially near their first order accumulation points. We first give definitions of standardized sets.

\begin{defn}[Standardized Sets]
Let $\Ii\subset\mathbb R$ be a set and $\gamma_0$ a real number. We say that $\Ii$ is \emph{standardized near} $\gamma_0$ if there exist a positive real number $\epsilon$, a positive integer $m$, and real numbers $b_1,\dots,b_m$, such that
$$\Ii\cap (\gamma_0-\epsilon,\gamma_0+\epsilon)\subset\left\{\gamma_0+\frac{b_i}{n}\mid i,n\in\mathbb N^+,1\leq i\leq m\right\}.$$
We say that $\Ii$ is 
\begin{enumerate}
    \item \emph{weakly standardized} if $\Ii$ is standardized near any $\gamma_0\in\bar\Ii\backslash\partial^2\Ii$,  and
    \item \emph{standardized} if $\partial^k\Ii$ is weakly standardized for any non-negative integer $k$ and $\partial^l\Ii=\emptyset$ for some positive integer $l$.
\end{enumerate}
\end{defn}
By Lemma \ref{lem: union of standardized sets} below, $\Ii$ is weakly standardized (resp. standardized) if and only if $\partial^0\Ii=\bar\Ii$ is weakly standardized (resp. standardized). 

Roughly speaking, a standardized set $\Ii$ has a filtration
$$\bar\Ii=\partial^0\Ii\supset\partial\Ii\supset\partial^2\Ii\dots\supset\partial^l\Ii=\emptyset$$
which consists of weakly standardized sets. Note that when $\Ii\subset [0,1]$ and $1$ is the only possible accumulation point of $\Ii$, a standardized set is always a subset of a hyperstandard set (\cite[3.2]{PS09}). This is the reason why we adopt the word ``standardized".

The main theorem of our paper is as follows.

\begin{thm}[\textbf{Main Theorem}]\label{thm: standardized lcts}
For any positive integer $d$, the set $\lct(d)$ is standardized.
\end{thm}

We remark that Theorem \ref{thm: standardized lcts} also holds for pairs which allow more complicated boundary coefficients. See Theorem \ref{thm: standardization lcts with boundary} for more details.

To obtain a better understanding of Theorem \ref{thm: standardized lcts}, we provide several examples below.

\begin{ex}\label{ex: lct of fermat hypersurfaces}
Let $d,a_1,\dots,a_d$ be fixed positive integers. For any positive integer $n$, we consider the ``diagonal" polynomial $f_n=x_1^{a_1}+x_2^{a_2}+\dots+x_{d}^{a_{d}}+x_{d+1}^{n}$ and the divisor $S_n:=(f_n=0)$ on $X:=\mathbb C^{d+1}$. Suppose that $\sum_{i=1}^d\frac{1}{a_i}<1$ and $n\gg 0$, then it is well known that $$\gamma_n:=\lct(X,0;S_n)=\min\left\{1,\sum_{i=1}^d\frac{1}{a_i}+\frac{1}{n}\right\}=\sum_{i=1}^d\frac{1}{a_i}+\frac{1}{n}.$$
Let $\gamma_0:=\sum_{i=1}^d\frac{1}{a_i}$. Then $\gamma_0$ is the accumulation point of $\{\gamma_n\}_{n=1}^{+\infty}$. Moreover,  $\gamma_n$ approaches $\gamma_0$ in a ``standardized way" as
$$\gamma_n=\gamma_0+\frac{1}{n}.$$
 In other words, $\{\gamma_n\}_{n=1}^{+\infty}$ is standardized near $\gamma_0$, hence $\{\gamma_n\}_{n=1}^{+\infty}$ is a standardized set. In fact, it is not hard to check that the set of lcts of ``diagonal" polynomials of dimension $d+1$,
$$\left\{\sum_{i=1}^{d+1}\frac{1}{c_i}\mid c_1,\dots,c_{d+1}\in\mathbb N^+\right\}\cap [0,1],$$
is a standardized set.
\end{ex}

\begin{ex}\label{ex: standardized ht2}
It is known that the set of lcts on $\mathbb C^2$ is
$$\mathcal{HT}_2=\left\{\frac{c_1+c_2}{c_1c_2+a_1c_2+a_2c_1}\mid a_1+c_1\geq\max\{2,a_2\},a_2+c_2\geq\max\{2,a_1\},a_1,a_2,c_1,c_2\in\mathbb N\right\}.$$
and the set of lcts on $\mathbb C^1$ is
$$\mathcal{HT}_1=\lct(1)=\left\{\frac{1}{k}\mid k\in\mathbb N^+\right\}\cup\{0\}$$
(cf. \cite[(15.5)]{Kol08}).
We have $\partial\mathcal{HT}_2=\mathcal{HT}_1\backslash\{1\}=\left\{\frac{1}{k}\mid k\in\mathbb N^+,k\geq 2\right\}\cup\{0\}$
and $\partial^2\mathcal{HT}_2=\{0\}$ (cf. \cite[Theoerem 7]{Kol08}). It is not hard to see that, for any $k\in\mathbb N^+$ and any sequence $\{\gamma_n\}_{n=1}^{+\infty}\subset\mathcal{HT}_2$ such that $$0<\frac{1}{k}=\gamma_0:=\lim_{n\rightarrow+\infty}\gamma_n,$$
possibly by passing to a subsequence and switching $c_1$ and $c_2$, we have
$$\gamma_n=\frac{c_{1,n}+c_2}{c_{1,n}c_2+a_1c_2+a_2c_{1,n}}$$
for some fixed $a_1,a_2,c_2$ and strictly increasing sequence of integers $c_{1,n}$, such that $a_2+c_2=k$ and $a_1\leq k$. Therefore, $\gamma_n$ approaches $\gamma_0$ in a ``standardized way" as
$$\gamma_n=\gamma_0+\frac{c_2(k-a_1)}{k}\cdot\frac{1}{kc_{1,n}+a_1c_2}\in\left\{\gamma_0+\frac{c_2(k-a_1)}{m}\mid m\in\mathbb N^+\right\}.$$
It is not hard to deduce that $\mathcal{HT}_2$ is standardized near any $\gamma_0\in\partial\mathcal{HT}_2\backslash\partial^2\mathcal{HT}_2$.

On the other hand, it is clear that the values in $\mathcal{HT}_2$ may not approach $0$ in a standardized way, hence $\mathcal{HT}_2$ is not standardized near $0$. Nevertheless, $0$ is a second order accumulation point of $\mathcal{HT}_2$, hence $\mathcal{HT}_2$ is still a weakly standardized set. Moreover, since $0$ is an accumulation point of $\mathcal{HT}_1$ and $\mathcal{HT}_1$ is standardized near $0$, we know that $\mathcal{HT}_2$ is standardized.
\end{ex}

Theorem \ref{thm: standardized lcts} could be potentially applied to the study on Han's uniform boundedness conjecture of minimal log discrepancies (cf. \cite[Conjecture 7.2]{HLL22}), especially on its weaker version for fixed germs cf. \cite[Conjecture 1.1]{MN18}, as the accumulation points of lcts naturally appear in the study of these conjectures. 

We also expect Theorem \ref{thm: standardized lcts} to be useful when estimating the precise values of log canonical thresholds in high dimension, especially the $1$-gap of lc thresholds.

Nevertheless, with Theorem \ref{thm: standardized lcts} settled, it will be interesting to ask whether other invariants in birational geometry behave similarly, such as the minimal log discrepancy and the canonical threshold. We will confirm that the sets of these invariants are standardized in some special cases in the following.

\medskip

\noindent\textbf{Minimal log discrepancies}. The minimal log discrepancy (mld for short) is another fundamental invariant in algebraic geometry, which measures how singular a variety is. The smaller the mld is, the worse the singularity is.

\begin{defn}
Let $(X\ni x)$ be an lc singularity. We define
$$\mld(X\ni x):=\min\{a(E,X)\mid E\text{ is over }X\ni x\}$$
to be the mld of $(X\ni x)$. The set of mlds for varieties of dimension $d$ is defined as
$$\mld(d):=\{\mld(X\ni x)\mid \dim X=d, X\ni x\text{ is lc}\}.$$
\end{defn}

Similar to the sets of log canonical thresholds, the sets of minimal log discrepancies of fixed dimension are also conjectured to satisfy the ACC \cite[Problem 5]{Sho88}, and its accumulation points are expected to come from lower dimension (cf. \cite[Version 1, Remark 1.2]{HLS19}). Unfortunately, the ACC conjecture for mlds is open in dimension $\geq 3$, hence it will be difficult to show the standardized behavior of this invariant. Nevertheless, we are able to prove that the sets of mlds are standardized in some special cases.

\begin{thm}\label{thm: standardized mlds}
\begin{enumerate}
    \item[$\mathrm{(1)}$] The sets $\mld(1)$ and $\mld(2)$ are standardized.
    \item[$\mathrm{(2)}$] The set $\{\mld(X)\mid\dim X=3,X\text{ is canonical}\}$ is standardized.
    \item[$\mathrm{(3)}$] For any positive integer $d$ and positive real number $\epsilon$, $$\{\mld(X\ni x)\mid \dim X=d, X\ni x\text{ has an }\epsilon\text{-plt blow up}\}$$
    is standardized. In particular,
    $$\{\mld(X\ni x)\mid \dim X=d, X\ni x\text{ is exceptional}\}$$
    is standardized.
\end{enumerate}
\end{thm}

We remark that Theorem \ref{thm: standardized mlds} also holds for lc pairs whose boundary coefficients belong to a finite set. See Section 3 for more details.

We also remark that Theorem \ref{thm: standardized mlds}(3) is actually important in the proof of Theorem \ref{thm: standardized lcts}. Note that the standardized behavior of mlds is very important due to the following example, which shows that a conjectural standardized behavior is already very helpful in the study of the ACC conjecture for mlds. In fact, the standardized behavior of lcts and mlds was observed when the first author examined the following example in \cite{LL22}.

\begin{ex}\label{ex: standardized 5/6}
A recent work of the first author and Luo shows that $\frac{5}{6}$ is the second largest accumulation point of global mlds in dimension $3$ \cite[Theorem 1.3]{LL22}. An important ingredient of the proof is \cite[Theorem 3.5]{LL22}, which essentially uses the conjectural standardized behavior of $\mld(3)$ in $\left(\frac{5}{6},1\right)$. The idea is as follows. 

For any real number $a\in\left(\frac{5}{6},1\right)$, we associate infinite equations to $a$ such that $a\in\mld(3)$ (almost) only if these equations have a common solution (cf. \cite[Definition 3.4]{LL22}). Denote the set of these equations by $\mathcal{E}(a)$. For each fixed $a$, the equations in $\mathcal{E}(a)$ are computable. As there are infinitely many equations, in practice, it is not hard to verify that the equations in $\mathcal{E}(a)$ do not have any common solution by checking finitely many of them, but it is difficult to verify that the equations in $\mathcal{E}(a)$ have a common solution. Moreover, since there are uncountably many real numbers $a$ in $\left(\frac{5}{6},1\right)$, we cannot consider all $\mathcal{E}(a)$ at the same time. To resolve these issues, a key idea is to decompose $\left(\frac{5}{6},1\right)$ as a disjoint union of subsets
$$\left(\frac{5}{6},1\right)=\cup_{n}\Ii_n\cup\tilde\Ii$$
such that 
\begin{enumerate}
    \item $\tilde\Ii$ satisfies the ACC and only accumulates to $\frac{5}{6}$ (hence these values will not influence the proof of \cite[Theorem 3.5]{LL22}), and
    \item for any fixed $n$, $\Ii_n$ is an open interval, and the equations in $\cap_{a\in\Ii_n}\mathcal{E}(a)$ do not have any common solution. This implies that $\Ii_n\cap\mld(3)=\emptyset$.
\end{enumerate}
The difficulty is that we need to guess what $\tilde\Ii$ is. Since the equations in $\mathcal{E}(a)$ have a common solution (almost) whenever $a\in\mld(3)\cap \left(\frac{5}{6},1\right)$ and $\mld(3)\cap \left(\frac{5}{6},1\right)$ is known to be an infinite set, we need to find a regular pattern of the values in $\tilde\Ii$.

At this point, a key observation is that the equations in $\mathcal{E}(a)$ heavily rely on the denominator of $a$. This is the reason why we conjectured that the denominator of $a$ grows standardly respect to $a-\frac{5}{6}$ when $a$ approaches $\frac{5}{6}$. With this conjecture in mind, an attempt on setting $$\tilde\Ii=\left\{\frac{5n+m}{6n+m}\mid m,n\in\mathbb N^+,1\leq m\leq 5\right\}\cup\left\{\frac{12}{13}\right\}$$
is successful. In summary, the conjecture on the standardized behavior of the set of mlds was essentially applied in the proof of \cite[Theorem 3.5]{LL22}, in a way that we directly ``guess the set of mlds out". Note that $\tilde\Ii$ is standardized as the only accumulation point of $\tilde\Ii$ is $\frac{5}{6}$ and
$$\frac{5n+m}{6n+m}=\frac{5}{6}+\frac{m}{36n+6m}\in\left\{\frac{5}{6}+\frac{m}{l}\mid l,m\in\mathbb N^+, 1\leq m\leq 5\right\}.$$

Similar strategies in \cite{LL22} can be applied to further studies on $\mld(3)$, especially for those values that are $\geq\frac{1}{2}$ as $\frac{1}{2}$ is the conjectured largest second order accumulation point of threefold mlds. See \cite{Rei87,Jia21,LX21} for related results.
\end{ex}

\noindent\textbf{Canonical thresholds}. The canonical threshold is another important invariant in birational geometry. In particular, canonical thresholds in dimension $3$ are deeply related to Sarkisov links in dimension $3$ (cf. \cite{Cor95,Pro18}). It is known that, in dimension $\leq 3$, the set of canonical thresholds satisfies the ACC \cite{Che22a,Che22b,HLL22}, and its accumulation points come from lower dimension \cite{Che22b,HLL22}. We show that the set of canonical thresholds is also standardized in dimension $\leq 3$, which actually follows from the proofs in \cite{Che22b,HLL22}:

\begin{thm}\label{thm: standardized threefold thresholds}
The set of canonical thresholds $\ct(d)$ in dimension $d$ is standardized when $d\leq 3$.
\end{thm}

\medskip

\noindent\textbf{Further discussions}. Yuchen Liu informed us that an invariant in K-stability and wall-crossing theory, the \emph{K-semistable threshold (walls)}, may also behave in a standardized way.

\begin{ex}[{\cite[Theorem 5.16]{ADL23}}]
The list of K-moduli walls of $\overline{\mathfrak{M}}^K_c$ (the K-moduli stack which parametrizes K-polystable log Fano pairs $(X,cD)$ admitting a $\Qq$-Gorenstein smoothing to $(\mathbb P^3,cS)$ where $S$ is a quartic surface) is
$$\left\{1-\frac{4}{n}\mid n\in\{6,8,10,12,13,14,16,18,22\}\right\},$$
which is a subset of a hyperstandard set (Definition \ref{defn: DCC and ACC}). Although this is a very specific example and the value of the walls are finite, using a hyperstandard set to describe the thresholds is natural in this case, and it is possible that larger classes of K-semistable thresholds also behave in a standardized way.
\end{ex}

\medskip

\noindent\textbf{Acknowledgement}. The authors would like to thank Guodu Chen, Christopher D. Hacon, Jingjun Han, Junpeng Jiao, Yuchen Liu, and Yujie Luo for useful discussions. The third author was partially supported by the NSF research grants no: DMS1952522, DMS-1801851 and by a grant from the Simons Foundation; Award Number: 256202.

\section{Preliminaries}

We adopt the standard notation and definitions in \cite{KM98,BCHM10} and will freely use them. 

\subsection{Sets}

\begin{defn}\label{defn: DCC and ACC}
Let $\Ii\subset\Rr$ be a set. We say that
\begin{enumerate}
    \item $\Ii$ satisfies the \emph{descending chain condition} (DCC) if any decreasing sequence in $\Ii$ stabilizes,
    \item $\Ii$ satisfies the \emph{ascending chain condition} (ACC) if any increasing sequence in $\Ii$ stabilizes,
    \item $\Ii$ is \emph{the standard set} if $\Ii=\{1-\frac{1}{n}\mid n\in\mathbb N^+\}\cup\{1\}$, and
    \item (\cite[3.2]{PS09}, \cite[2.2]{Bir19}) $\Ii$ is a \emph{hyperstandard set} if there exists a finite set $\Ii_0\subset\mathbb R_{\geq 0}$ such that $0,1\in\Ii_0$ and $\Ii=\{1-\frac{\gamma}{n}\mid n\in\mathbb N^+,\gamma\in\Ii_0\}\cap [0,1]$.
\end{enumerate}
\end{defn}

\subsection{Pairs and singularities}

\begin{defn}[Pairs, {cf. \cite[Definition 3.2]{CH21}}] \label{defn sing}
A \emph{pair} $(X/Z\ni z, B)$ consists of a contraction $\pi: X\rightarrow Z$, a (not necessarily closed) point $z\in Z$, and an $\mathbb{R}$-divisor $B\geq 0$ on $X$, such that $K_X+B$ is $\Rr$-Cartier over a neighborhood of $z$. If $\pi$ is the identity map and $z=x$, then we may use $(X\ni x, B)$ instead of $(X/Z\ni z,B)$. In addition, if $B=0$, then we use $X\ni x$ instead of $(X\ni x,0)$. If $(X\ni x,B)$ is a pair for any codimension $\geq 1$ point $x\in X$, then we call $(X,B)$ a pair. A pair $(X\ni x, B)$ is called a \emph{germ} if $x$ is a closed point.
\end{defn}

\begin{defn}[Singularities of Pairs]\label{defn: relative mld}
 Let $(X\ni x,B)$ be a pair and $E$ a prime divisor over $X$ such that $x\in \Center_XE$. Let $f: Y\rightarrow X$ be a log resolution of $(X,B)$ such that $\Center_Y E$ is a divisor, and suppose that $K_Y+B_Y=f^*(K_X+B)$ over a neighborhood of $x$. We define $a(E,X,B):=1-\mult_EB_Y$ to be the \emph{log discrepancy} of $E$ with respect to $(X,B)$.
 
 For any prime divisor $E$ over $X$, we say that $E$ is \emph{over} $X\ni x$ if $\Center_XE=\bar x$. We define
 $$\mld(X\ni x,B):=\inf\{a(E,X,B)\mid E\text{ is over }X\ni x\}$$
 to be the \emph{minimal log discrepancy} (\emph{mld}) of $(X\ni x,B)$. We define $$\mld(X,B):=\inf\{a(E,X,B)\mid E\text{ is exceptional over }X\}.$$
 We define
 $$\tmld(X,B):=\inf\{a(E,X,B)\mid E\text{ is over }X\}$$
to be the \emph{total minimal log discrepancy (tmld)} of $(X,B)$.
 
 Let $\epsilon$ be a non-negative real number. We say that $(X\ni x,B)$ is \emph{lc} (resp. \emph{klt, $\epsilon$-lc, $\epsilon$-klt}) if $\mld(X\ni x,B)\geq 0$ (resp. $>0$, $\geq\epsilon$, $>\epsilon$). We say that $(X,B)$ is \emph{lc} (resp. \emph{klt, $\epsilon$-lc, $\epsilon$-klt}) if $\tmld(X,B)\geq 0$ (resp. $>0$, $\geq\epsilon$, $>\epsilon$). 
 
 We say that $(X,B)$ is \emph{canonical} (resp. \emph{terminal}, \emph{plt}, \emph{$\epsilon$-plt}) if $\mld(X,B)\geq 1$ (resp. $>1,>0,>\epsilon$).
 \end{defn}

\begin{defn}\label{defn: alct local}  Let $a$ be a non-negative real number, $(X\ni x,B)$ (resp. $(X,B)$) an lc pair, and $D\geq 0$ an $\Rr$-Cartier $\Rr$-divisor on $X$. We define
$$a\text{-}\lct(X\ni x,B;D):=\sup\{-\infty,t\mid t\geq 0, (X\ni x,B+tD)\text{ is }a\text{-lc}\}$$
$$\text{(resp. }a\text{-}\lct(X,B;D):=\sup\{-\infty,t\mid t\geq 0, (X,B+tD)\text{ is }a\text{-lc}\} \text{)}$$
to be the \emph{$a$-lc threshold} of $D$ with respect to $(X\ni x,B)$ (resp. $(X,B)$). We define
$$\ct(X\ni x,B;D):=\sup\{-\infty,t\mid t\geq 0, (X\ni x,B+tD)\text{ is }1\text{-lc}\}$$
$$\text{(resp. }\ct(X,B;D):=\sup\{-\infty,t\mid t\geq 0, (X,B+tD)\text{ is canonical}\} \text{)}$$
to be the \emph{canonical threshold} of $D$ with respect to $(X\ni x,B)$ (resp. $(X,B)$). We define $\lct(X\ni x,B;D):=0\text{-}\lct(X\ni x,B;D)$ (resp. $\lct(X,B;D):=0\text{-}\lct(X,B;D)$) to be the \emph{lc threshold} of $D$ with respect to $(X\ni x,B)$ (resp. $(X,B)$).
\end{defn}

\begin{defn}
Assume that $X$ is a normal variety and $B$ is an $\Rr$-divisor on $X$. We write $B\in\Ii$ if the coefficients of $B$ belong to $\Ii$. For any positive integer $d$, we define
$$\mld(d,\Ii):=\{\mld(X\ni x,B)\mid \dim X=d,(X\ni x,B)\text{ is lc}, B\in\Ii\},$$
$$\lct(d,\Ii):=\{\lct(X,B;D)\mid \dim X=d,(X,B)\text{ is lc}, B\in\Ii,D\in\mathbb N^+\},$$
and
$$\ct(d,\Ii):=\{\ct(X,B;D)\mid \dim X=d,(X,B)\text{ is canonical}, B\in\Ii,D\in\mathbb N^+\}.$$

We let $\mld(0,\Ii)=\lct(0,\Ii)=\ct(0,\Ii):=\{0\}$. For any non-negative integer $d$, we let $\mld(d):=\mld(d,\{0\})$, $\lct(d):=\lct(d,\{0,1\})$, and $\ct(d):=\ct(d,\{0,1\})$. 
\end{defn}

 \subsection{Complements}\label{section: complements}
\begin{defn}\label{defn: complement}
Let $n$ be a positive integer, $\Ii_0\subset (0,1]$ a finite set, and $(X/Z\ni z,B)$ and $(X/Z\ni z,B^+)$ two pairs. We say that $(X/Z\ni z,B^+)$ is an \emph{$\Rr$-complement} of $(X/Z\ni z,B)$ if 
\begin{itemize}
    \item $(X/Z\ni z,B^+)$ is lc,
    \item $B^+\geq B$, and
    \item $K_X+B^+\sim_{\Rr}0$ over a neighborhood of $z$.
\end{itemize}
We say that $(X/Z\ni z,B^+)$ is an \emph{$n$-complement} of $(X/Z\ni z,B)$ if
\begin{itemize}
\item $(X/Z\ni z,B^+)$ is lc,
\item $nB^+\geq \lfloor (n+1)\{B\}\rfloor+n\lfloor B\rfloor$, and
\item $n(K_X+B^+)\sim 0$ over a neighborhood of $z$.
\end{itemize}
 We say that $(X/Z\ni z,B)$ is $\Rr$-complementary if $(X/Z\ni z,B)$ has an $\Rr$-complement. We say that $(X/Z\ni z,B^+)$ is a \emph{monotonic $n$-complement} of $(X/Z\ni z,B)$ if $(X/Z\ni z,B^+)$ is an $n$-complement of $(X/Z\ni z,B)$ and $B^+\geq B$.
 
 We say that $(X/Z\ni z,B^+)$ is an \emph{$(n,\Ii_0)$-decomposable $\Rr$-complement} of $(X/Z\ni z,B)$ if there exist a positive integer $k$, $a_1,\dots,a_k\in\Ii_0$, and  $\Qq$-divisors $B_1^+,\dots,B_k^+$ on $X$, such that
\begin{itemize}
\item $\sum_{i=1}^ka_i=1$ and  $\sum_{i=1}^ka_iB_i^+=B^+$,
\item $(X/Z\ni z,B^+)$ is an $\Rr$-complement of $(X/Z\ni z,B)$, and
\item  $(X/Z\ni z,B_i^+)$ is an $n$-complement of itself for each $i$.
\end{itemize}
\end{defn}

\begin{thm}[{\cite[Theorem 1.10]{HLS19}}]\label{thm: ni decomposable complement}
Let $d$ be a positive integer and $\Ii\subset [0,1]$ a DCC set. Then there exist a positive integer $n$ and a finite set $\Ii_0\subset (0,1]$ depending only on $d$ and $\Ii$ satisfying the following. 

Assume that $(X/Z\ni z,B)$ is a pair of dimension $d$ and $B\in\Ii$, such that $X$ is of Fano type over $Z$ and $(X/Z\ni z,B)$ is $\Rr$-complementary. Then $(X/Z\ni z,B)$ has an $(n,\Ii_0)$-decomposable $\Rr$-complement. Moreover, if $\bar\Ii\subset\Qq$, then $(X/Z\ni z,B)$ has a monotonic $n$-complement.
\end{thm}

\subsection{Plt blow-ups}
\begin{defn}\label{defn: reduced component}
	Let $(X\ni x,B)$ be a klt germ and $\epsilon$ a positive real number. A \emph{plt} (resp. \emph{$\epsilon$-plt}) \emph{blow-up} of $(X\ni x,B)$ is a divisorial contraction $f: Y\rightarrow X$ with a prime exceptional divisor $E$ over $X\ni x$, such that $(Y/X\ni x,f^{-1}_*B+E)$ is plt (resp. $\epsilon$-plt) and $-E$ is ample over $X$.
\end{defn}

\begin{lem}[{\cite[3.1]{Sho96},\cite[Proposition 2.9]{Pro00},\cite[Theorem 1.5]{Kud01},\cite[Lemma 1]{Xu14}}]\label{lem: existence plt blow up}
	Assume that $(X\ni x,B)$ is a klt germ such that $\dim X\geq 2$. Then there exists a plt blow-up of $(X\ni x,B)$.
\end{lem}

\begin{defn}
Let $(X\ni x,B)$ be an lc germ. We say that $(X\ni x,B)$ is \emph{exceptional} if for any $\mathbb R$-divisor $G\geq 0$ on $X$ such that $(X\ni x,B+G)$ is lc, there exists at most one lc place of $(X\ni x,B+G)$.
\end{defn}

\subsection{Special sets}

\begin{defn}
Let $\Ii\subset [0,1]$ be a set, $d$ a positive integer, and $c$ a positive real number. We define
$$\Ii_+:=\left(\{0\}\cup\left\{\sum_{i=1}^n\gamma_i\mid \gamma_1,\dots,\gamma_n\in\Ii\right\}\right)\cap [0,1],$$
$$D(\Ii):=\left\{\frac{m-1+\gamma}{m}\mid m\in\mathbb N^+,\gamma\in\Ii_+\right\},$$
$$D(\Ii,c):=\left\{\frac{m-1+\gamma+kc}{m}\mid m,k\in\mathbb N^+,\gamma\in\Ii_+\right\}\cap [0,1],$$
\begin{align*}
    \mathfrak{N}(d,\Ii,c):=\Bigg\{(X,B)\Biggm|
    \begin{array}{r@{}l}
        (X,B)\text{ is projective lc}, K_X+B\equiv 0,\dim X=d,\\
     B=L+C,L\in D(\Ii),0\not=C\in D(\Ii,c)
    \end{array}\Bigg\},
    \end{align*}
    $$\mathfrak{R}(d,\Ii,c):=\left\{(X,B)\mid (X,B)\in\mathfrak{N}(n,\Ii,c), (X,B) \text{ is }\Qq\text{-factorial klt}, \rho(X)=1,1\leq n\leq d\right\},$$
    $$N(d,\Ii):=\{c\mid c\in [0,1],\mathfrak{N}(d,\Ii,c)\not=\emptyset\},$$
    and
    $$K(d,\Ii):=\{c\mid c\in [0,1],\mathfrak{R}(d,\Ii,c)\not=\emptyset\}.$$
    We define $N(0,\Ii)=K(0,\Ii):=\{\frac{1-\gamma}{n}\mid \gamma\in\Ii_+, n\in\mathbb N^+\}\cup\{0\}$ and $N(-1,\Ii)=K(-1,\Ii):=\{0\}$.
\end{defn}
The following results in \cite{HMX14} are used in the proof of the main theorem.

\begin{thm}[{\cite[Lemmas 11.2, 11.4, Proposition 11.5, Theorem 1.11]{HMX14}}]\label{thm: hmx accumulation point preliminaries}
Let $d$ be a non-negative integer and $\Ii\subset [0,1]$ a set. Then
\begin{enumerate}
    \item[$\mathrm{(1)}$] $\lct(d,\Ii)\subset\lct(d+1,\Ii)$ and $N(d-1,\Ii)\subset N(d,\Ii)$,
    \item[$\mathrm{(2)}$] $N(d,\Ii\cup\{1\})=K(d,\Ii)$,
    \item[$\mathrm{(3)}$] if $\Ii=\Ii_+$, then $\lct(d,\Ii)=N(d-1,\Ii)$, and
    \item[$\mathrm{(4)}$] if $1\in\Ii,\Ii=\Ii_+$, and $\partial\Ii\subset\{1\}$, then $\partial\lct(d+1,\Ii)=\lct(d,\Ii)\backslash\{1\}$.
\end{enumerate}
\end{thm}

\subsection{Basic properties of standardized sets}

The behavior of standardized sets is generally similar to the behavior of DCC sets and ACC sets. However, there are still some differences. 

\begin{ex}
It is clear that any subset of a DCC (resp. ACC) set is still DCC (resp. ACC). However, a subset of a standardized set may no longer be standardized. Consider the sets $\Ii_1:=\{\frac{n}{n^2+1}\mid n\in\mathbb N^+\}$ and $\Ii_2=\Ii_1\cup\{\frac{1}{n}+\frac{1}{m}\mid n,m\in\mathbb N^+\}$. Then $\Ii_1\subset\Ii_2$. It is not hard to check that $\Ii_2$ is a standardized set but $\Ii_1$ is not. This is because although $\Ii_1$ and $\Ii_2$ are both not standardized near $0$, $0\not\in\partial^2\Ii_1$, but $0\in\partial^2\Ii_2$.
\end{ex}

We summarize the following properties on standardized sets below which we will use in this paper.

\begin{lem}\label{lem: basic properties standardized}
Let $\Ii$ be a set of real numbers and $\gamma_0$ a real number. Then:
\begin{enumerate}
    \item[$\mathrm{(1)}$] If $\gamma_0\not\in\partial\Ii$, then $\Ii$ is standardized near $\gamma_0$.
    \item[$\mathrm{(2)}$] If $\gamma_0\in\partial^2\Ii$, then $\Ii$ is not standardized near $\gamma_0$.
    \item[$\mathrm{(3)}$] For any real number $a$, $\Ii$ is standardized near $\gamma_0$ if and only if $\{\gamma+a\mid \gamma\in\Ii\}$ is standardized near $a+\gamma_0$.
    \item[$\mathrm{(4)}$] For any non-zero number $c$, $\Ii$ is standardized near $\gamma_0$ if and only if $\{c\gamma\mid \gamma\in\Ii\}$ is standardized near $c\gamma_0$.
    \item[$\mathrm{(5)}$] Suppose that $\Ii=\cup_{i=1}^k\Ii_i$. Then $\Ii$ is standardized near $\gamma_0$ if and only if each $\Ii_i$ is standardized near $\gamma_0$. In particular, $\Ii$ is standardized near $\gamma_0$ if and only if any subset of $\Ii$ is standardized near $\gamma_0$.
    \item[$\mathrm{(6)}$] $\Ii$ is standardized near $\gamma_0$ if and only if $\Ii\cap (\gamma_0-\epsilon_0,\gamma_0+\epsilon_0)$ is standardized near $\gamma_0$ for some positive real number $\epsilon_0$.
    \item[$\mathrm{(7)}$] If $\{\gamma-\gamma_0\mid\gamma\in\Ii\}\subset\mathbb Q$, then $\Ii$ is standardized near $\gamma_0$ if and only if there exists a positive integer $I$ and a positive real number $\epsilon$, such that $$\Ii\cap(\gamma_0-\epsilon,\gamma_0+\epsilon)\subset\left\{\gamma_0+\frac{I}{n}\mid n\in\mathbb Z\backslash\{0\}\right\}\cup\{\gamma_0\}.$$
    \item[$\mathrm{(8)}$] $\Ii$ is standardized near $\gamma_0$ if and only if $\bar\Ii$ is standardized near $\gamma_0$.
\end{enumerate}
\end{lem}
\begin{proof}
For any $\gamma_0\not\in\partial\Ii$, we may pick a positive real number $\epsilon$ such that $(\gamma_0-\epsilon,\gamma_0+\epsilon)\cap\Ii=\{\gamma_0\}$ or $\emptyset$. This implies (1). 

Suppose that $\Ii$ is standardized near $\gamma_0$ for some $\gamma_0\in\partial^2\Ii$. Then there exist a positive real number $\epsilon$, a positive integer $m$, and real numbers $b_1,\dots,b_m$, such that $\Ii\cap (\gamma_0-\epsilon,\gamma_0+\epsilon)\subset\left\{\gamma_0+\frac{b_i}{n}\mid i,n\in\mathbb N^+,1\leq i\leq m\right\}$. Thus the only accumulation point of $\Ii\cap (\gamma_0-\epsilon,\gamma_0+\epsilon)$ is $\gamma_0$, hence $\gamma_0\not\in\partial^2\Ii$, a contradiction. This implies (2).

(3)(4)(6) are obvious.

We prove (5). If $\Ii$ is standardized near $\gamma_0$, then there exist a positive real number $\epsilon$, a positive integer $m$, and non-zero real numbers $b_1,\dots,b_m$, such that
$$\Ii_i\cap (\gamma_0-\epsilon,\gamma_0+\epsilon)\subset\Ii\cap (\gamma_0-\epsilon,\gamma_0+\epsilon)\subset\left\{\gamma_0+\frac{b_j}{n}\mid j,n\in\mathbb N^+,1\leq j\leq m\right\}.$$
for each $1\leq i\leq k$, hence $\Ii_i$ is standardized near $\gamma_0$ for each $i$. If $\Ii_i$ is standardized near $\gamma_0$ for each $i$, then there exist positive integers $m_1,\dots,m_k$, a finite set of real numbers $\{b_{i,j}\}_{1\leq i\leq k,1\leq j\leq m_i}$, and real numbers $\epsilon_1,\dots,\epsilon_k$, such that
$$\Ii_i\cap (\gamma_0-\epsilon_i,\gamma_0+\epsilon_i)\subset\left\{\gamma_0+\frac{b_{i,j}}{n}\mid j,n\in\mathbb N^+,1\leq j\leq m_i\right\}$$
for each $i$. Let $\epsilon:=\min\{\epsilon_1,\dots,\epsilon_k\}$, then
$$\Ii\cap (\gamma_0-\epsilon,\gamma_0+\epsilon)\subset\left\{\gamma_0+\frac{b_{i,j}}{n}\mid i,j,n\in\mathbb N^+,1\leq i\leq k,1\leq j\leq m_i\right\},$$
hence $\Ii$ is standardized.

We prove (7). The if part is obvious. Suppose that $\Ii$ is standardized near $\gamma_0$, then there exist a positive real number $\epsilon$, a positive integer $m$, and non-zero real numbers $b_1,\dots,b_m$, such that $$\Ii\cap(\gamma_0-\epsilon,\gamma_0+\epsilon)\subset\left\{\gamma_0+\frac{b_i}{n}\mid i,n\in\mathbb N^+, 1\leq i\leq m\right\}.$$
Since $\{\gamma-\gamma_0\mid \gamma\in\Ii\}\subset\mathbb Q$, we may assume that $b_i\in\mathbb Q$ for each $i$. We may let $I$ be a common denominator of the elements in $\left\{\frac{1}{b_i}\mid b_i\not=0\right\}$.

The if part of (8) follows from (5). Suppose that $\Ii$ is standardized near $\gamma_0$, then there exist a positive real number $\epsilon$, a positive integer $m$, and non-zero real numbers $b_1,\dots,b_m$,  such that $$\Ii\cap(\gamma_0-\epsilon,\gamma_0+\epsilon)\subset\left\{\gamma_0+\frac{b_i}{n}\mid i,n\in\mathbb N^+, 1\leq i\leq m\right\}.$$
Thus 
$$\bar\Ii\cap (\gamma_0-\epsilon,\gamma_0+\epsilon)\subset\overline{\Ii\cap (\gamma_0-\epsilon,\gamma_0+\epsilon)}\subset\left\{\gamma_0+\frac{b_i}{m}\mid i,n\in\mathbb N^+,1\leq i\leq n\right\}\cup\{\gamma_0\},$$
hence $\bar\Ii$ is standardized near $\gamma_0$.
\end{proof}

\begin{lem}\label{lem: standardized iff subsequence standardized}
Let $\Ii$ be a set of real numbers and $\gamma_0$ a real number. Suppose that for any sequence $\{\gamma_i\}_{i=1}^{+\infty}\subset\Ii$ such that $\lim_{i\rightarrow+\infty}\gamma_i=\gamma_0$, there exists an infinite subsequence of  $\{\gamma_i\}_{i=1}^{+\infty}$ which is standardized near $\gamma_0$. Then $\Ii$ is standardized near $\gamma_0$.
\end{lem}
\begin{proof}
For any non-zero real number $\gamma$, we let $[\gamma]$ be its $\Qq$-class under multiplication: $[\gamma]=[\gamma']$ if and only if $\gamma=s\gamma'$ for some $s\in\mathbb Q^{\times}$. By Lemma \ref{lem: basic properties standardized}(3), possibly by replacing $\gamma_0$ with $0$ and $\Ii$ with $\{\gamma-\gamma_0\mid\gamma\in\Ii\}$, we may assume that $\gamma_0=0$. For any positive real number $\epsilon$, we consider $\Ii_{\mathbb Q,\epsilon}:=\{[\gamma]\mid \gamma\in (\Ii\cap  (-\epsilon,\epsilon))\backslash\{0\}\}$.

Suppose that $\Ii_{\mathbb Q,\epsilon}$ is an infinite set for any positive real number $\epsilon$. Then we may pick a sequence $\{\gamma_i\}_{i=1}^{+\infty}\subset\Ii$ such that $[\gamma_i]\not=[\gamma_j]$ for any $i\not=j$ and $\lim_{i\rightarrow+\infty}\gamma_i=0$. By assumption, there exists a strictly increasing sequence of integers $\{r_i\}_{i=1}^{+\infty}$ such that $\{\gamma_{r_i}\}_{i=1}^{+\infty}$ is standardized near $0$, so there exist a positive integer $m$ and real numbers $b_1,\dots,b_m$, such that
$$\{\gamma_{r_i}\}_{i=1}^{+\infty}\subset\left\{\frac{b_j}{n}\mid j,n\in\mathbb N^+,1\leq j\leq m\right\}.$$
This is not possible as the $\Qq$-classes of  $\{\frac{b_j}{n}\mid j,n\in\mathbb N^+,1\leq j\leq m\}$ are finite. Therefore, $\Ii_{\mathbb Q,\epsilon_0}$ is a finite set for some positive real number $\epsilon_0$. By Lemma \ref{lem: basic properties standardized}(6), we may replace $\Ii$ with $\Ii\cap (-\epsilon_0,\epsilon_0)$ and write $\Ii=\cup_{i=1}^k\Ii_i$ for some positive integer $k$, such that $[\gamma_i]\not=[\gamma_j]$ for any $\gamma_i\in\Ii_i$ and $\gamma_j\in\Ii_j$ and any $i\not=j$, and $[\alpha]=[\beta]$ for any $\alpha,\beta\in\Ii_i$ and any $i$.  By Lemma \ref{lem: basic properties standardized}(5), we only need to show that $\Ii_i$ is standardized near $0$ for any $i$. Therefore, we may assume that $k=1$. In particular, there exists a non-zero real number $c$ and a set $\Ii'\subset\mathbb Q$ such that $\Ii=\{c\gamma'\mid\gamma'\in\Ii'\}$. By Lemma \ref{lem: basic properties standardized}(4), possibly by replacing $\Ii$ with $\Ii'$, we may assume that $\Ii\subset\mathbb Q$.

Suppose that $\Ii$ is not standardized near $0$. By Lemma \ref{lem: basic properties standardized}(7), there exists a sequence $\{\frac{p_i}{q_i}\}_{i=1}^{+\infty}\subset\Ii$ such that $\gcd(p_i,q_i)=1$, $\lim_{i\rightarrow+\infty}|p_i|=+\infty$, and $\lim_{i\rightarrow+\infty}\frac{p_i}{q_i}=0$. It is clear that no infinite subsequence of $\{\frac{p_i}{q_i}\}_{i=1}^{+\infty}$ is standardized near $0$, a contradiction.
\end{proof}
\begin{lem}\label{lem: union of standardized sets}
Let $\Ii,\Ii'$ be two sets of real numbers. Then:
\begin{enumerate}
    \item $\Ii$ is weakly standardized (resp. standardized) if and only if $\bar\Ii$ is weakly standardized (resp. standardized).
    \item If $\Ii$ and $\Ii'$ are weakly standardized (resp. standardized), then $\Ii\cup\Ii'$ is weakly standardized (resp. standardized).
\end{enumerate}
\end{lem}
\begin{proof}
Since $\bar\Ii\backslash\partial^2\Ii=\bar\Ii\backslash\partial^2\bar\Ii=\bar{\bar\Ii}\backslash\partial^2\bar\Ii$ and $\partial^k\Ii=\partial^k\bar\Ii$ for any non-negative integer $k$, (1) follows from Lemma \ref{lem: basic properties standardized}(8). 

Since $\partial^k(\Ii\cup\Ii')=\partial^k\Ii\cup\partial^k\Ii'$ for any non-negative integer $k$, (2) follows from Lemma \ref{lem: basic properties standardized}(5).
\end{proof}

We summarize some additional properties of standardized sets in the following lemma. The lemma is interesting in its own right. However, we do not need this lemma in the rest of this paper.

\begin{lem}\label{lem: basic property standardized sets}
Let $\Ii,\Ii'$ be two sets of real numbers.
\begin{enumerate}
    \item[$\mathrm{(1)}$] $\Ii$ is weakly standardized if and only if $\Ii$ is standardized near any $\gamma_0\in\partial\Ii\backslash\partial^2\Ii$.
    \item[$\mathrm{(2)}$] If $\Ii,\Ii'$ are standardized and both satisfy the DCC (resp. ACC), then $$\Ii'':=\{\gamma+\gamma'\mid \gamma\in\Ii,\gamma'\in\Ii'\}$$ is standardized.
    \item[$\mathrm{(3)}$] If $\Ii'$ is a finite set, then $\Ii$ is standardized if and only if $\Ii\cup\Ii'$ is standardized.
    \item[$\mathrm{(4)}$] If $\Ii'$ is an interval, $\Ii$ is DCC or ACC, and $\Ii$ is standardized, then $\Ii\cap\Ii'$ is standardized.
    \item[$\mathrm{(5)}$] If $\Ii\subset [0,1]$ satisfies the DCC and is standardized, then $\Ii_+$ and $D(\Ii)$ are standardized.
    \item[$\mathrm{(6)}$] If $\Ii\subset [0,+\infty)$ satisfies the ACC and is standardized, then $\{\frac{\gamma}{n}\mid\gamma\in\Ii,n\in\mathbb N^+\}$ is standardized.
\end{enumerate}
\end{lem}
\begin{proof}
(1) It follows from Lemma \ref{lem: basic properties standardized}(1).

(2) Since $\Ii$ and $\Ii'$ both satisfy the DCC (resp. ACC), we have $$\partial^k\Ii''=\cup_{i=0}^k\{\gamma+\gamma'\mid \gamma\in\partial^i\Ii,\gamma'\in\partial^{k-i}\Ii'\}.$$
Since $\Ii$ and $\Ii'$ are standardized, $\partial^{l}\Ii=\emptyset$ and $\partial^{l'}\Ii'=\emptyset$ for some positive integers $l,l'$. Thus $\partial^{l+l'}\Ii''=\emptyset$. By induction on $l+l'$ and Lemma \ref{lem: union of standardized sets}(2), we only need to show that $\Ii''$ is weakly standardized. By Lemma \ref{lem: standardized iff subsequence standardized} and (1), we only need to show that for any $\gamma_0''\in\partial\Ii''\backslash\partial^2\Ii''$ and any sequence $\{\gamma_i''\}_{i=1}^{+\infty}$ such that $\lim_{i\rightarrow+\infty}\gamma_i''=\gamma_0''$, a subsequence of $\{\gamma_i''\}_{i=1}^{+\infty}$ is standardized near $\gamma_0''$. We may write $\gamma_i''=\gamma_i+\gamma_i'$ where $\gamma_i\in\Ii$ and $\gamma_i'\in\Ii'$. Possibly by passing to a subsequence, we may assume that $\gamma_i,\gamma_i'$ are increasing (resp. decreasing), $\lim_{i\rightarrow+\infty}\gamma_i=\gamma_0$, and $\lim_{i\rightarrow+\infty}\gamma_i'=\gamma_0'$. If $\{\gamma_i\}_{i=1}^{+\infty}$ and $\{\gamma_i'\}_{i=1}^{+\infty}$ both have strictly increasing (resp. strictly decreasing) subsequences, then possibly by passing to subsequences, we may assume that $\{\gamma_i\}_{i=1}^{+\infty}$ and $\{\gamma_i'\}_{i=1}^{+\infty}$ are strictly increasing (resp. strictly decreasing). Since
$$\gamma_0''=\lim_{j\rightarrow+\infty}(\gamma_0+\gamma_j')=\lim_{j\rightarrow+\infty}\lim_{i\rightarrow+\infty}(\gamma_i+\gamma_j'),$$ 
 $\gamma_0''\in\partial^2\Ii''$, a contradiction. Thus possibly by passing to a subsequence and switching $\Ii,\Ii'$, we may assume that $\gamma_i=\gamma_0$ for each $i$. By Lemma \ref{lem: basic properties standardized}(3), $\{\gamma_i''\}_{i=1}^{+\infty}$ is standardized near $\gamma_0''$, and we get (2). 

(3) We have $\partial^k(\Ii\cup\Ii')=\partial^k\Ii$ for any positive integer $k$. For any real number $\gamma_0$ and non-negative integer $k$, since $\Ii'$ is a finite set, by Lemma \ref{lem: basic properties standardized}(5), $\partial^k\Ii$ is standardized near $\gamma_0$ if and only if $\partial^k(\Ii\cup\Ii')$ is standardized near $\gamma_0$. This implies (3).

(4) By Lemma \ref{lem: basic properties standardized}(4), possibly by replacing $\Ii$ with $\{-\gamma\mid\gamma\in\Ii\}$ and $\Ii'$ with $\{-\gamma'\mid\gamma'\in\Ii'\}$, we may assume that $\Ii$ is DCC. Let $a:=\inf\Ii'\in\{-\infty\}\cup\mathbb R$ and $c:=\sup\Ii'\in\{+\infty\}\cup\mathbb R$.  Since $\Ii$ satisfies the DCC, $a$ is not an accumulation point of $\Ii\cap\Ii'$. Thus $\partial^k(\Ii\cap\Ii')=(\partial^k\Ii\cap\Ii')\backslash\{a,c\}$ or $(\partial^k\Ii\cap\Ii'\backslash\{a\})\cup\{c\}$ for any positive integer $k$. By (3) and induction on the minimal non-negative integer $l$ such that $\partial^l\Ii=\emptyset$, we only need to show that $\Ii\cap\Ii'$ is weakly standardized. By (1) and Lemma \ref{lem: basic properties standardized}(6), we only need to show that $\Ii\cap\Ii'$ is standardized near $c$ when $c<+\infty$ and $c\in\partial(\Ii\cap\Ii')\backslash\partial^2(\Ii\cap\Ii')$. Since $\Ii$ satisfies the DCC, $c\in\partial(\Ii\cap\Ii')\backslash\partial^2(\Ii\cap\Ii')$ if and only if $c\in\partial\Ii\backslash\partial^2\Ii$, hence $\Ii$ is standardized near $c$ when $c\in\partial(\Ii\cap\Ii')\backslash\partial^2(\Ii\cap\Ii')$. Statement (4) follows from Lemma \ref{lem: basic properties standardized}(5).

(5) Suppose that $\Ii\subset [0,1]$ satisfies the DCC and is standardized. First we show that $\Ii_+$ is standardized. We let $\Ii_1:=\Ii$ and let $\Ii_k:=\{\gamma+\tilde\gamma\mid \gamma\in\Ii,\tilde\gamma\in\Ii_{k-1}\}$ for any integer $k\geq 2$. Since $\Ii$ satisfies the DCC, we may let $\bar\gamma:=\min\{1,\gamma\in\Ii\mid \gamma>0\}$. By (2), $\Ii_k$ is standardized for any positive integer $k$. Since $\Ii_+=(\Ii_{\lfloor\frac{1}{\bar\gamma}\rfloor}\cup\{0\})\cap [0,1]$, by (3)(4), $\Ii_+$ is standardized.

Now we show that $D(\Ii)$ is standardized. We may replace $\Ii$ with $\Ii_+$ and suppose that $\Ii=\Ii_+$. Then $(\partial^k\Ii)_+=\partial^k\Ii\cup\{0\}$ for any non-negative integer $k$. We have
$$D(\Ii)=\left\{\frac{m-1+\gamma}{m}\mid m\in\mathbb N^+,\gamma\in\Ii\right\}.$$
Let $k_0$ be the minimal positive integer such that $\partial^{k_0}\Ii\not=\emptyset$. By induction, we have
$$\partial^k D(\Ii)=\{1\}\cup\{\frac{m-1+\gamma}{m}\mid m\in\mathbb N^+,\gamma\in\partial^k\Ii\}$$ 
for any $1\leq k\leq k_0$, $\partial^{k_0+1}D(\Ii)=\{1\}$, and $\partial^{k}D(\Ii)=\emptyset$ for any $k\geq k_0+2$. By (3) and induction on $k_0$, we only need to show that $D(\Ii)$ is weakly standardized. There are two cases:

\medskip

\noindent\textbf{Case 1}. The set $\Ii$ is a finite set. Then $1$ is the only accumulation point of $D(\Ii)$, and it is clear that $D(\Ii)$ is standardized near $1$. By (1), $D(\Ii)$ is weakly standardized.

\medskip

\noindent\textbf{Case 2}. The set $\Ii$ is not a finite set. Then $1\in\partial^2D(\Ii)$. For any $c\in [0,1)$ and any sequence $\{c_i\}_{i=1}^{+\infty}\subset D(\Ii)$ such that $\lim_{i\rightarrow+\infty}c_i=c$, possibly by passing to a subsequence we have
$$c_i=\frac{m-1+\gamma_i}{m}$$
such that $m$ is a constant and $\gamma_i\in\Ii$ for each $i$. If $c\not\in\partial^2D(\Ii)$, then by Lemma \ref{lem: basic properties standardized}(3)(4), $\{c_i\}_{i=1}^{+\infty}$ is standardized near $c$. By Lemma \ref{lem: standardized iff subsequence standardized}, $D(\Ii)$ is standardized near $c$, hence $D(\Ii)$ is weakly standardized.

(6) Since $\Ii$ satisfies the ACC, $\Ii\subset [0,M]$ for some positive integer $M$. By Lemma \ref{lem: basic properties standardized}(4), possibly by replacing $\Ii$ with $\{\frac{\gamma}{M}\mid \gamma\in\Ii\}$, we may assume that $M=1$. Let $\Ii':=\{1-\gamma\mid\gamma\in\Ii\}$, then $\Ii'\subset [0,1]$ satisfies the DCC and is standardized. By the same argument as in (5) and Lemma \ref{lem: basic properties standardized}(3)(4),
$$\left\{\frac{\gamma}{n}\mid n\in\mathbb N^+,\gamma\in\Ii\right\}=\left\{1-\frac{m-1+\gamma'}{m}\mid m\in\mathbb N^+,\gamma'\in\Ii'\right\}$$
is standardized.
\end{proof}

\section{Standardization of (some) log discrepancies}

\begin{thm}\label{thm: finite coefficient surface mld standardized}
Let $\Ii\subset [0,1]$ be a finite set. Then $\mld(1,\Ii)$ and $\mld(2,\Ii)$ are standardized.
\end{thm}
\begin{proof}
Since $\mld(1,\Ii)=\{1,1-\gamma\mid\gamma\in\Ii\}$ is a finite set, $\mld(1,\Ii)$ is standardized.

We first show that $\mld(2,\Ii)$ is standardized near any $a_0>0$. Fix $a_0>0$ and let $\{a_i\}_{i=1}^{+\infty}\subset\mld(2,\Ii)$ be any sequence such that $\lim_{i\rightarrow+\infty}a_i=a_0$.  Let $(X_i\ni x_i,B_i)$ be a surface singularity such that $B_i\in\Ii$ and $\mld(X_i\ni x_i,B_i)=a_i$. 

\begin{claim}\label{claim: abd chl21}
Possibly by passing to a subsequence, we have
$$a_i=\frac{\alpha A_i+\beta}{A_i+\delta}$$
where $\alpha,\delta\geq 0$ and $\beta>0$ are constants such that $\delta\in\mathbb Q$, and $A_i\in\mathbb N^+$.  
\end{claim}
\begin{proof}
    This essentially follows from \cite[Lemma A.2]{CH21} (see also \cite[Lemma 3.3]{Ale93}), but since the argument of \cite[Lemma A.2]{CH21} is very long, we provide a short proof here. 

    Let $\epsilon:=\min\{a_0,\Ii_{>0}\}$. If the possibilities of the dual graphs of the minimal resolution of $X_i\in x_i$ is finite, then there are only finitely many possibilities of $a_i=\mld(X_i\in x_i,B_i)$, which is not possible. Therefore, by \cite[Lemma A.6]{CH21}, possibly by passing to a subsequence, we may assume that 
$$a_i=\mld(X_i\ni x_i,B_i)=\pld(X_i\ni x_i,B_i).$$ 
Possibly by passing to a subsequence, we may assume that $(X_i\ni x_i,B_i)$ is of one of the types as in \cite[Lemma A.2(1)]{CH21} ($\mathfrak{F}_{\epsilon,\Ii}$), \cite[Lemma A.2(2)]{CH21} ($\mathfrak{C}_{\epsilon,\Ii}$), or \cite[Lemma A.2(3)]{CH21} ($\mathfrak{T}_{\epsilon,\Ii}$). If $(X_i\ni x_i,B_i)$ is of type $\mathfrak{F}_{\epsilon,\Ii}$ for each $i$, then there are finitely many possibilities of the dual graph of the minimal resolution of $X_i\ni x_i$, which is again not possible. If $(X_i\ni x_i,B_i)$ is of type $\mathfrak{T}_{\epsilon,\Ii}$ for each $i$, then \cite[Lemma A.2(3)]{CH21} implies that
$$a_i=\pld(X_i\ni x_i,B_i)=\frac{\alpha_i}{m_i-q_i}$$
where $q_i<m_i\leq\lfloor\frac{2}{\epsilon}\rfloor^{\lfloor\frac{2}{\epsilon}\rfloor}$, $q_i$ and $m_i$ are integers, and $\alpha_i$ belongs to a finite set as it is constructed as in \cite[Page 35, Line 17]{CH21} and $\Ii$ is a finite set. In this case, $\frac{\alpha_i}{m_i-q_i}$ belongs to a finite set, which is not possible. Therefore, we may assume that $(X_i\ni x_i,B_i)$ is of type $\mathfrak{C}_{\epsilon,\Ii}$ for each $i$. \cite[Lemma A.2(2)]{CH21} implies that
$$a_i=\pld(X_i\ni x_i,B_i)=\frac{(A_i+\frac{m_{2,i}}{m_{2,i}-q_{2,i}})\frac{\alpha_{1,i}}{m_{1,i}-q_{1,i}}+\frac{q_{1,i}}{m_{1,i}-q_{1,i}}\cdot \frac{\alpha_{2,i}}{m_{2,i}-q_{2,i}}}{A_i+\frac{q_{1,i}}{m_{1,i}-q_{1,i}}+\frac{m_{2,i}}{m_{2,i}-q_{2,i}}}$$
where $q_{1,i}<m_{1,i}\leq\lfloor\frac{2}{\epsilon}\rfloor^{\lfloor\frac{2}{\epsilon}\rfloor}$ and $q_{2,i}<m_{2,i}\leq\lfloor\frac{2}{\epsilon}\rfloor^{\lfloor\frac{2}{\epsilon}\rfloor}$, $q_{1,i},q_{2,i},m_{1,i},m_{2,i}$ are integers, and $A_i$ is a non-negative integers. Therefore, possibly passing to a subsequence, we may assume that $q_{1,i},q_{2,i},m_{1,i},m_{2,i}$ are constants, and $A_i>0$. We may let
$$\alpha:=\frac{\alpha_{1,i}}{m_{1,i}-q_{1,i}},\beta:=\frac{m_{2,i}}{m_{2,i}-q_{2,i}}\cdot\frac{\alpha_{1,i}}{m_{1,i}-q_{1,i}}+\frac{q_{1,i}}{m_{1,i}-q_{1,i}}\cdot \frac{\alpha_{2,i}}{m_{2,i}-q_{2,i}},$$
and
$$\delta:=\frac{q_{1,i}}{m_{1,i}-q_{1,i}}+\frac{m_{2,i}}{m_{2,i}-q_{2,i}}.$$
\end{proof}
\noindent\textit{Proof of Theorem \ref{thm: finite coefficient surface mld standardized}} continued. Let $\alpha,\beta,\delta$ and $A_i$ be as in Claim \ref{claim: abd chl21}. Then $\alpha=a_0$, and
$$a_i=a_0+\frac{\beta-a_0\delta}{A_i+\delta}.$$
Therefore, $\{a_i\}_{i=1}^{+\infty}$ is standardized near $a_0$. By Lemma \ref{lem: standardized iff subsequence standardized}, $\mld(2,\Ii)$ is standardized near $a_0$.

By  \cite[Lemma A.2]{CH21} (see also \cite[Lemma 3.3]{Ale93}, we have
$$\left\{0,\frac{1}{n}\mid n\in\mathbb N^+\right\}\subset\partial\mld(2,\Ii)\subset\left\{0,\frac{1-\gamma}{n}\mid n\in\mathbb N^+, \gamma\in\Ii_+\right\}$$
and
$$\partial^2\mld(2,\Ii)=\{0\}.$$
Since $\Ii$ is a finite set, $\Ii_+$ is a finite set, hence $\{0,\frac{1-\gamma}{n}\mid\gamma\in\Ii_+,n\in\mathbb N^+\}$ is standardized near $0$. By Lemma \ref{lem: basic properties standardized}(5), $\partial\mld(2,\Ii)$ is standardized near $0$. Since $\overline{\mld(2,\Ii)}\subset [0,+\infty)$, $\mld(2,\Ii)$ is standardized.
\end{proof}

\begin{thm}\label{thm: finite coefficient terminal threefold mld standardized}
Let $\Ii\subset [0,1]$ be a finite set. Then
$$\Ii':=\{\mld(X,B)\mid \dim X=3,B\in\Ii\}\cap [1,+\infty)$$
is standardized and its only accumulation point is $1$.
\end{thm}
\begin{proof}
By \cite[Corollary 1.5]{Nak16}, $1$ is the only accumulation point of $\Ii'$, so we only need to show that $\Ii'$ is standardized near $1$. By Lemma \ref{lem: standardized iff subsequence standardized}, we only need to show that for any sequence of pairs $\{(X_i,B_i)\}_{i=1}^{+\infty}$ such that $\dim X_i=3,B_i\in\Ii$, and $\mld(X_i,B_i)\geq 1$, $\{\mld(X_i,B_i)\}_{i=1}^{+\infty}$ has a subsequence which is standardized near $1$. Possibly by passing to a subsequence and replacing each $(X_i,B_i)$ with a $\mathbb Q$-factorialization, we may assume that each $X_i$ is $\mathbb Q$-factorial. By \cite[Theorem 6.12]{HLL22}, possibly by passing to a subsequence, we may find a positive integer $l$ depending only on $\Ii$, and prime divisors $E_i$ that are exceptional over $X_i$, such that $a(E_i,X_i,B_i)=\mld(X_i,B_i)>1$ and $a(E_i,X_i,0)\leq 1+\frac{l}{I_i}$, where $I_i$ is the Cartier index of $K_{X_i}$ near the generic point $x_i$ of $\Center_{X_i}E_i$. By \cite[Corollary 5.2]{Kaw88}, for any prime divisor $D$ on $X_i$ , $I_iD$ is Cartier near $x_i$. Since $\Ii$ is a finite set, possibly by passing to a subsequence, we may assume that $a(E_i,X_i,B_i)=1+\frac{\gamma}{I_i}$ where $\gamma\in (0,l]$ is a constant. It is clear that $\{a(E_i,X_i,B_i)\}_{i=1}^{+\infty}$ is standardized near $1$ and the theorem follows.
\end{proof}

\begin{lem}\label{lem: complicated plt blow up}
Let $d$ be a positive integer, $\epsilon$ and $c$ two positive real numbers, and $\Ii\subset [0,1]$ and $\Ii'\subset [0,+\infty)\cap\mathbb Q$ two finite sets. Then there exist a finite set $\Ii_1\subset (0,+\infty)$ and a finite set $\Ii_2\subset [0,+\infty)\cap\mathbb Q$ depending only on $d,\epsilon,c,\Ii$, and $\Ii'$ which satisfy the following. 

Assume that $(X\ni x,B)$ is an lc pair of dimension $d$, such that
\begin{enumerate}
    \item $B=\Delta+sS$ such that $\Delta\in\Ii$, $S\in\Ii'$, and
    \item $(X\ni x,\Delta+cS)$  has an $\epsilon$-plt blow-up $f: Y\rightarrow X$ which extracts a prime divisor $E$.
\end{enumerate}
Then 
$a(E,X,B)=\frac{\alpha-(s-c)\beta}{n},$
where $\alpha\in\Ii_1,\beta\in\Ii_2$, and $n\in\mathbb N^+$.
\end{lem}
\begin{proof}
Possibly by replacing $c,s$ and $\Ii'$, we may assume that $S$ is a Weil divisor. By cutting $X$ by general hyperplane sections and applying induction on dimension, we may assume that $x$ is a closed point. Let $\Delta_Y,B_Y$, and $S_Y$ be the strict transforms of $\Delta,B$ and $S$ on $Y$ respectively, and $a:=a(E,X,B)$.
$$K_Y+B_Y+(1-a)E=f^*(K_X+B).$$
Since $(X\ni x,B)$ is lc, $a\geq 0$. Let
$$K_E+B_E:=(K_Y+B_Y+E)|_E\quad\text{and}\quad K_E+B_E':=(K_Y+\Delta_Y+cS_Y+E)|_E.$$
Since $a\geq 0$ and $-E$ is ample$/X$, $-(K_E+B_E)$ is nef. Since $f$ is an $\epsilon$-plt blow-up of $(X\ni x,\Delta+cS)$, $(E,B_E')$ is an $\epsilon$-klt log Fano pair. By \cite[Theorem 1.1]{Bir21}, $E$ belongs to a bounded family. Thus there exist a positive integer $M$ depending only on $d$ and $\epsilon$, and a very ample divisor $H$ on $E$, such that $-K_E\cdot H^{d-2}\leq M$. By adjunction (cf. \cite[Theorem 3.10]{HLS19}), we may write
$$B_E=\sum_D\frac{m_D-1+\gamma_D+sk_D}{m_D}D\quad\text{and}\quad B'_E=\sum_D\frac{m_D-1+\gamma_D+ck_D}{m_D}D$$
where the sums are taken over all prime divisors $D$ on $E$, $m_D$ are positive integers, $k_D$ are non-negative integers, and $\gamma_D\in\Ii_+$. Since $(E,B_E')$ is $\epsilon$-klt and $\Ii\subset [0,1]$ is a finite set, $\gamma_D$ belongs to a finite set of non-negative real numbers, $m_D$ belongs to a finite set of positive integers, and $k_D$ belongs to a finite set of non-negative integers.

Since $0<-(K_E+B'_E)\cdot H^{d-2}\leq M$ and $D\cdot H^{d-2}$ is a positive integer for each $D$, $-(K_E+B_E)\cdot H^{d-2}$ is of the form $\alpha'-(s-c)\beta'$, where $\alpha'=-(K_E+B_E')\cdot H^{d-2}$ belongs to a finite set of positive real numbers and $\beta':=(\sum_D\frac{k_D}{m_D}D)\cdot H^{d-2}$ belongs to a finite set of non-negative rational numbers.

Let $H_1,\ldots,H_{d-2}$ be general elements in $|H|$, $C:=E\cap H_1\cap H_2\cdots \cap H_{d-2}$, and $r:=\lfloor\frac{1}{\epsilon}\rfloor!$. Since $(Y/X\ni x,B_Y+E)$ is $\epsilon$-plt, by \cite[Theorem 3.10]{HLS19}, $rE|_{E}$ is a $\Qq$-Cartier Weil divisor. In particular, $-E\cdot C$ belongs to the discrete set $\frac{1}{r}\mathbb N^+$. Since $(K_Y+B_Y+(1-a)E)\cdot C=0$, we have
$$a=\frac{-(K_Y+B_Y+E)\cdot C}{-E\cdot C}=\frac{-(K_E+B_E)\cdot H^{d-2}}{-E\cdot C}.$$
Thus $a=\frac{r(\alpha'-(s-c)\beta')}{n}$, where $n\in\mathbb N^+$. We may let $\alpha:=r\alpha'$ and $\beta:=r\beta'$.
\end{proof}

\begin{lem}\label{lem: eplt kc ld standardized}
Let $d$ be a positive integer, $\epsilon$ a positive real number, and $\Ii\subset [0,1]$ a finite set. Then
\begin{align*}
  \Bigg\{a(E,X,B)\Biggm|
    \begin{array}{r@{}l}
        (X\ni x,B)\text{ has an }\epsilon\text{-plt blow-up } f: Y\rightarrow X \\
    \text{ which extracts }E,\dim X=d,B\in\Ii
    \end{array}\Bigg\}
    \end{align*}
    is standardized and its only possible accumulation point is $0$.
\end{lem}
\begin{proof}
It follows from Lemma \ref{lem: complicated plt blow up} by letting $c=s=1$ and $S=0$.
\end{proof}

\begin{thm}\label{thm: finite coefficient eplt mld standardized}
Let $d$ be a positive integer, $\epsilon$ a positive real number, and $\Ii\subset [0,1]$ a finite set. Then
$$\Ii_1(d,\epsilon,\Ii):=\{\mld(X\ni x,B))\mid \dim X=d, (X\ni x,B)\text{ has an }\epsilon\text{-plt blow up}\}$$
    is standardized and its only possible accumulation point is $0$. In particular,
    $$\Ii_2(d,\Ii):=\{\mld(X\ni x,B)\mid \dim X=d, (X\ni x,B)\text{ is exceptional}\}$$
    is standardized and its only possible accumulation point is $0$.
\end{thm}
\begin{proof}
 By \cite[Theorems 1.2, 1.3]{HLS19}, the only possible accumulation point of $\Ii_1(d,\epsilon,\Ii)$ and $\Ii_2(d,\Ii)$ is $0$. By Lemma \ref{lem: basic properties standardized}(1), we only need to show that $\Ii_1(d,\epsilon,\Ii)$ and $\Ii_2(d,\Ii)$ are standardized near $0$. 
 
 By \cite[Lemma 3.22]{HLS19}, there exists a positive real number $\epsilon'$ depending only on $d$ and $\Ii$, such that for any exceptional pair $(X\ni x,B)$ of dimension $d$ with $B\in\Ii$, $(X\ni x,B)$ has an $\epsilon'$-plt blow-up. In particular, $\Ii_2(d,\Ii)\subset\Ii_1(d,\epsilon',\Ii)$. 
 
 By \cite[Theorem 1.3]{HLS19} and Lemma \ref{lem: eplt kc ld standardized}, $\Ii_1(d,\epsilon,\Ii)\cap [0,\epsilon]$ and $\Ii_1(d,\epsilon',\Ii)\cap [0,\epsilon']$ are standardized near $0$. By Lemma \ref{lem: basic properties standardized}(6), $\Ii_1(d,\epsilon,\Ii)$ and $\Ii_1(d,\epsilon',\Ii)$ are standardized near $0$. By Lemma \ref{lem: basic properties standardized}(5), $\Ii_2(d,\Ii)$ is standardized near $0$, and we are done.
\end{proof}

\begin{proof}[Proof of Theorem \ref{thm: standardized mlds}]
It follows from Theorems \ref{thm: finite coefficient surface mld standardized}, \ref{thm: finite coefficient terminal threefold mld standardized} and \ref{thm: finite coefficient eplt mld standardized} when $\Ii=\{0\}$.
\end{proof}

\section{Standardization of log canonical thresholds}

\begin{lem}\label{lem: coeff 1 imply low dimensional lct}
Let $d$ be a positive integer, $c$ a non-negative real number, and $\Ii\subset [0,1]$ a set such that $1\in\Ii_+$. Suppose that $(X,B)\in\mathfrak{N}(d,\Ii,c)$ is a pair such that $(X,B)$ is not klt. Then $c\in N(d-1,\Ii)$.
\end{lem}
\begin{proof}
Possibly by replacing $(X,B)$ with a dlt modification, we may assume that $\lfloor B\rfloor\not=0$ and $X$ is $\Qq$-factorial. We may assume that $c\not=0$. We may write $B=L+C$ such that $L\in D(\Ii)$ and $0\not=C\in D(\Ii,c)$. If $\lfloor C\rfloor\not=0$, then
$$\frac{m-1+\gamma+kc}{m}=1$$
for some $m,k\in\mathbb N^+$ and $\gamma\in\Ii_+$. By Theorem \ref{thm: hmx accumulation point preliminaries}(1), $c=\frac{1-\gamma}{k}\in N(0,\Ii)\subset N(d-1,\Ii)$. Thus we may assume that $\lfloor C\rfloor=0$.

If $\lfloor L\rfloor=0$, then
$$\frac{m_1-1+\gamma_1}{m_1}+\frac{m_2-1+\gamma_2+kc}{m_2}=1$$
for some $m_1,m_2,k\in\mathbb N^+$ and $\gamma_1,\gamma_2\in\Ii_+$. Since $c\not=0$, either $m_1=1$ or $m_2=1$. If $m_1=1$, then
$$c=\frac{1-\gamma_2-m_2\gamma_1}{k}\in N(0,\Ii)\subset N(d-1,\Ii),$$
and if $m_2=1$, then
$$c=\frac{1-\gamma_1-m_1\gamma_2}{m_1k}\in N(0,\Ii)\subset N(d-1,\Ii).$$
Thus we may assume that $\lfloor L\rfloor\not=0$. We let $T$ be an irreducible component of $\lfloor L\rfloor$. We run a $(K_X+L)$-MMP $\phi: X\dashrightarrow X'$ which terminates with a Mori fiber space $X'\rightarrow Z$. Then this MMP is $C$-positive, hence $C$ is not contracted by this MMP.

If $T$ is contracted by $\phi$, then there exists a step of the MMP $\psi: X''\rightarrow X'''$ which is a divisorial contraction and contracts the strict transform of $T$ on $X''$. Let $B'',L'',C'',T''$ be the strict transforms of $B,L,C,T$ on $X''$ respectively. Since $\psi$ is $C''$-positive, $T''$ intersects $C''$. Since $(X,L)$ is dlt, $(X'',L'')$ is dlt, hence $T''$ is normal. Let
$$K_{T''}+B_{T''}:=(K_{X''}+B'')|_{T''},$$
then $(T'',B_{T''})\in\mathfrak{N}(d-1,\Ii,c)$. Thus $c\in N(d-1,\Ii)$. Therefore, we may assume that $T$ is not contracted by $\phi$.

We let $B',L',C',T'$ be the strict transforms of $B,L,C,T$ on $X'$ respectively. Note that $T'$ is normal as $(X',L')$ is dlt. Since $\phi$ is $C$-positive, $C'$ dominates $Z$.  

If $\dim Z>0$, then we let $F$ be a very general fiber of $X'\rightarrow Z$, and let
$$K_F+B_F:=(K_{X'}+B')|_F,$$
then $(F,B_F)\in\mathfrak{N}(d-\dim Z,\Ii,c)$.  Thus $c\in N(d-\dim Z,\Ii)\subset N(d-1,\Ii)$. Thus we may assume that $\dim Z=0$ and $\rho(X')=1$.

If $d\geq 2$, then $T'$ intersects $C'$. Let
$$K_{T'}+B_{T'}:=(K_{X'}+B')|_{T'},$$
then $(T',B_{T'})\in\mathfrak{N}(d-1,\Ii,c)$. Thus $c\in N(d-1,\Ii)$ and we are done.

If $d=1$, then we have
$$\sum_{j=1}^{l_1}\frac{m_j-1+\gamma_j}{m_j}+\sum_{j=1}^{l_2}\frac{n_j-1+\gamma'_j+k_jc}{n_j}=1$$
for some $l_1\in\mathbb N$, $l_2,m_j,n_j,k_j\in\mathbb N^+$, and $\gamma_j,\gamma_j'\in\Ii_+$. Since $c>0$, possibly by reordering indices, either $m_j=1$ for every $j$ and $n_j=1$ for every $j\geq 2$, or $m_j=1$ for every $j\geq 2$ and $n_j=1$ for every $j$. Thus either
$$c=\frac{1-n_1(\sum_{j=1}^{l_1}\gamma_j+\sum_{j=2}^{l_2}\gamma_j')-\gamma_1'}{k_1+n_1\sum_{j=2}^{l_2}k_j}\in N(0,\Ii)\subset N(d-1,\Ii)$$
or
$$c=\frac{1-m_1(\sum_{j=2}^{l_1}\gamma_j+\sum_{j=1}^{l_2}\gamma_j')-\gamma_1}{m_1\sum_{j=1}^{l_2}k_j}\in N(0,\Ii)\subset N(d-1,\Ii)$$
and we are done.
\end{proof}

\begin{thm}\label{thm: standardization lcts with boundary}
Let $d$ be a non-negative integer and $\Ii\subset [0,1]$ a set, such that $1\in\Ii$, $\Ii=\Ii_+$, $1$ is the only possible accumulation point of $\Ii$, and $\Ii$ is standardized. Then $\lct(d,\Ii)$ is standardized.
\end{thm}
\begin{proof}

\textbf{Step 1}. In this step, we reduce our theorem to the case when $d\geq 2$ and show that we only need to prove that $\lct(d,\Ii)$ is standardized near any $c\in\lct(d-1,\Ii)\backslash(\lct(d-2,\Ii)\cup\{1\})$.

\smallskip

Since $1$ is the only possible accumulation point of $\Ii$ and $\Ii$ is standardized, there exists a positive integer $m$ and non-negative real numbers $b_1,\dots,b_m$, such that
$$\Ii\subset\left\{1-\frac{b_i}{n}\mid i,n\in\mathbb N^+, 1\leq i\leq m\right\}.$$
In particular, $\Ii$ satisfies the DCC. By \cite[Theorem 1.1]{HMX14}, $\lct(d,\Ii)$ satisfies the ACC for any non-negative integer $d$.

If $d=0$, then the theorem follows from the definition. If $d=1$, then 
\begin{align*}
 \lct(d,\Ii)&=\left\{\frac{1-\gamma}{k}\mid k\in\mathbb N^+, \gamma\in\Ii\right\}\subset\left\{\frac{b_i}{nk}\mid i,n,k\in\mathbb N^+, 1\leq i\leq m\right\}\\
 &\subset\left\{\frac{b_i}{n}\mid i,n\in\mathbb N^+,1\leq i\leq m\right\}. 
\end{align*}
Thus the only possible accumulation point of $\lct(d,\Ii)$ is $0$, and by Lemma \ref{lem: basic properties standardized}(5), $\lct(d,\Ii)$ is standardized near $0$. Thus $\lct(d,\Ii)$ is standardized and we are done. 

Therefore, we may assume that $d\geq 2$. By induction on dimension and Theorem \ref{thm: hmx accumulation point preliminaries}(4), we only need to show that $\lct(d,\Ii)$ is weakly standardized, that is, for any $c\in\partial\lct(d,\Ii)\backslash\partial^2\lct(d,\Ii)=\lct(d-1,\Ii)\backslash(\lct(d-2,\Ii)\cup\{1\})$, $\lct(d,\Ii)$ is standardized near $c$. 


\medskip

\noindent\textbf{Step 2}. For any $c\in\lct(d-1,\Ii)\backslash(\lct(d-2,\Ii)\cup\{1\})$ and $c_i\in\lct(d,\Ii)$ such that $\lim_{i\rightarrow+\infty}c_i=c$, we construct pairs $(X_i,B_i)\in \mathfrak{R}(d,\Ii,c_i)$ in this step.

\smallskip

For $c\in\lct(d-1,\Ii)\backslash(\lct(d-2,\Ii)\cup\{1\})$, we define
$$\epsilon_c:=\sup\{t\mid 0\leq t\leq 1, (c,c+t)\cap\lct(d-1,\Ii)=\emptyset\}.$$
Since $c\not\in\partial\lct(d-1,\Ii)$, $\epsilon_c>0$. By  Theorem \ref{thm: hmx accumulation point preliminaries}(4), $c$ is the only accumulation point of $(c,c+\epsilon_c)\cap\lct(d,\Ii)$. Since $\lct(d,\Ii)$ satisfies the ACC, by Lemma \ref{lem: basic properties standardized}(6), we only need to show that $(c,c+\epsilon_c)\cap\lct(d,\Ii)$ is standardized near $c$. By Lemma \ref{lem: standardized iff subsequence standardized}, we only need to show that for any sequence $\{c_i\}_{i=1}^{+\infty}\subset(c,c+\epsilon_c)\cap\lct(d,\Ii)$ such that $\lim_{i\rightarrow+\infty}c_i=c$, a subsequence of $\{c_i\}_{i=1}^{+\infty}$ is standardized near $c$. Possibly by passing to a subsequence, we may assume that $c_i$ is strictly decreasing.

In the following, we will fix $c\in\lct(d-1,\Ii)\backslash(\lct(d-2,\Ii)\cup\{1\})$ and a sequence $\{c_i\}_{i=1}^{+\infty}\subset (c,c+\epsilon_c)\cap\lct(d,\Ii)$ such that $\lim_{i\rightarrow+\infty}c_i=c$. In particular, $c_i\not\in\lct(d-1,\Ii)$ for each $i$. By Theorem \ref{thm: hmx accumulation point preliminaries}(2)(3), $c_i\in K(d-1,\Ii)\backslash K(d-2,\Ii)$ for each $i$ and $c\in K(d-2,\Ii)\backslash K(d-3,\Ii)$. Since $d\geq 2$ and $c\not\in\lct(d-2,\Ii)$, $c\not=0$. Thus there exists a sequence of pairs $(X_i,B_i=L_i+C_i)$, such that
\begin{enumerate}
    \item $\dim X_i=d-1$ and $\rho(X_i)=1$,
    \item $K_{X_i}+B_i\equiv 0$ and $(X_i,B_i)$ is $\Qq$-factorial klt, and
    \item $L_i\in D(\Ii)$ and $0\not=C_i\in D(\Ii,c_i)$.
\end{enumerate}
In particular, for each $i$, we may write
$$L_i=\sum_j\frac{m_{i,j}-1+\gamma_{i,j}}{m_{i,j}}L_{i,j}\quad\text{and}\quad C_i=\sum_{j}\frac{n_{i,j}-1+\gamma'_{i,j}+k_{i,j}c_i}{n_{i,j}}C_{i,j},$$
such that $m_{i,j},n_{i,j},k_{i,j}\in\mathbb N^+$, $\gamma_{i,j},\gamma_{i,j}'\in\Ii$, and $L_{i,j},C_{i,j}$ are prime divisors. We write $C_i=R_i+c_iS_i$ where $R_i:=\sum_{j}\frac{n_{i,j}-1+\gamma'_{i,j}}{n_{i,j}}C_{i,j}$ and $S_i:=\sum_{j}\frac{k_{i,j}}{n_{i,j}}C_{i,j}$. 

We let $a_i:=\tmld(X_i,L_i+R_i+cS_i)$. Possibly by passing to a subsequence, we may assume that $a_i$ is increasing or decreasing, and let $a:=\lim_{i\rightarrow+\infty}a_i$. 

\medskip

\noindent\textbf{Step 3}. In this step, we show that $a=0$.

\smallskip

Suppose that $a>0$. Then there exists a positive real number $\epsilon$ such that $(X_i,L_i+R_i+cS_i)$ is $\epsilon$-lc for each $i$, hence $X_i$ belongs to a bounded family by \cite[Theorem 1.1]{Bir21}. Thus there exist a positive integer $M$ which does not depend on $i$, and very ample divisors $H_i$ on $X_i$, such that $-K_{X_i}\cdot H_i^{d-2}\leq M$ for each $i$.

Since $(X_i,L_i+R_i+cS_i)$ is $\epsilon$-lc for each $i$, $\frac{m_{i,j}-1+\gamma_{i,j}}{m_{i,j}}\leq 1-\epsilon$ and $\frac{n_{i,j}-1+\gamma'_{i,j}+k_{i,j}c}{n_{i,j}}\leq 1-\epsilon$ for any $i,j$. Since $1$ is the only possible accumulation point of $\Ii$ and $c>0$, $m_{i,j},\gamma_{i,j},n_{i,j},\gamma_{i,j}',k_{i,j}$ belong to a finite set. Thus  $L_i\cdot H_i^{d-2}\geq 0,R_i\cdot H_i^{d-2}\geq 0$, and $S_i\cdot H_i^{d-2}>0$ belong to discrete sets. Since $K_{X_i}+B_i\equiv 0$, we have 
$$(K_{X_i}+L_i+R_i+c_iS_i)\cdot H_i^{d-2}=0.$$
Thus $c_i=\frac{p_i}{q_i}$, where $p_i=-(K_{X
_i}+L_i+R_i)\cdot H_i^{d-2}$ belongs to a finite set of positive real numbers, and $q_i=S_i\cdot H_i^{d-2}$ belongs to a discrete set of positive real numbers. Thus the only possible accumulation point of $\{c_i\}_{i=1}^{+\infty}$ is $0$, which is not possible as $c\not=0$.

Thus $a=0$. Let $a_i':=\tmld(X_i,B_i)$. Since $a=0$ and $0<a_i'\leq a_i$, $\lim_{i\rightarrow+\infty}a_i'=0$. Possibly by passing to a subsequence, we may assume that $a_i'$ is strictly decreasing.

\medskip

\noindent\textbf{Step 4}. In this step, we find a positive integer $N$, a finite set $\Ii_0\subset (0,1]$, a positive real number $\epsilon_0$, and divisors $T_i$ over $X_i$. We then construct $(N,\Ii_0)$-decomposable $\mathbb R$-complements $(X_i,L_i+R_i+cS_i+G_i)$ and Mori fiber spaces $(X_i',B_i')\rightarrow Z_i$, and reduce our theorem to the case when $\dim Z_i=0$ and $\rho(X_i')=1$.

\smallskip

By Theorem \ref{thm: ni decomposable complement}, there exist a positive integer $N$ and a finite set $\Ii_0\subset (0,1]$ depending only on $d,\Ii$ and $c$, such that for any $\mathbb R$-complementary pair $(X/Z\ni z,B)$ where $X$ is of Fano type over $Z$, $\dim X= d-1$, and $B\in D(\Ii\cup\{c\})$, $(X/Z\ni z,B)$ has an $(N,\Ii_0)$-decomposable $\Rr$-complement. We let $\epsilon_0:=\min\left\{\frac{\gamma_0}{2N}\mid \gamma_0\in\Ii_0\right\}>0$. Since $0=a=\lim_{i\rightarrow+\infty}a_i$, possibly by passing to a subsequence, we may assume that $a_i\leq\min\left\{\epsilon_0, 1\right\}$ for each $i$. We let $T_i$ be a prime divisor over $X_i$ such that $a(T_i,X_i,B_i)=\tmld(X_i,B_i)=a_i'$. We let $(X_i,L_i+R_i+cS_i+G_i)$ be an $(N,\Ii_0)$-decomposable $\Rr$-complement of $(X_i,L_i+R_i+cS_i)$. By our construction, $a(T_i,X_i,L_i+R_i+cS_i+G_i)=0$.

We construct a pair $(X_i',B_i')$ and a Mori fiber space $X_i'\rightarrow Z_i$ in the following way: \begin{itemize}
    \item If $T_i$ is on $X_i$, we let $Z_i:=\{pt\}$, $(X_i',B_i')=(X_i,B_i)$, $L_i':=L_i-L_i\wedge T_i$, $L'_{i,j}:=L_{i,j}$, $C'_i:=C_i-C_i\wedge T_i$, $C'_{i,j}:=C_{i,j}$, $R_i':=R_i-R_i\wedge T_i,S_i':=S_i-(\mult_{T_i}S_i)T_i$, and $G_i':=G_i-G_i\wedge T_i$. We let $T_i':=T_i$.
    \item If $T_i$ is exceptional over $X_i$, we let $f_i: Y_i\rightarrow X_i$ be a birational contraction which only extracts $T_i$, and let $B_{Y_i},L_{Y_i},L_{Y_i,j},C_{Y_i},C_{Y_i,j},R_{Y_i},S_{Y_i},G_{Y_i}$ be the strict transforms of $B_i,L_i,L_{i,j},C_i,C_{i,j},R_i,S_i,G_i$ on $Y_i$ respectively. We run a $(K_{Y_i}+B_{Y_i})$-MMP, which terminates with a Mori fiber space $X_i'\rightarrow Z_i$. We let $L_i',L'_{i,j},C'_i,C'_{i,j},R'_i,S'_i,G_i',T_i'$ be the strict transforms of $L_{Y_i},L_{Y_i,j},C_{Y_i},C_{Y_i,j},R_{Y_i},S_{Y_i},G_{Y_i},T_i$ on $X_i'$ respectively, and let $B_i':=(1-a'_i)T_i'+L_i'+C_i'$.
\end{itemize}
By our construction, $T_i'\not=0$, $T_i'$ dominates $Z_i$, $B_i'=(1-a'_i)T_i'+L_i'+C_i'$, and $C_i'=R_i'+c_iS_i'$. Moreover, since $K_{X_i}+L_i+R_i+cS_i+G_i\equiv 0$ and $a(T_i,X_i,L_i+R_i+cS_i+G_i)=0$, $(X_i,L_i+R_i+cS_i+G_i)$ and $(X_i',T_i'+L_i'+R_i'+cS_i'+G_i')$ are crepant. Since $K_{X_i}+L_i+R_i+c_iS_i\equiv 0$ and $a(T_i,X_i,L_i+R_i+c_iS_i)=a_i'$, $(X_i',(1-a_i')T'_i+L_i'+R'_i+c_iS'_i)$ and $(X_i,B_i)$ are crepant. Thus $K_{X'_i}+(1-a_i')T'_i+L_i'+R'_i+c_iS'_i\equiv 0$ and $(X_i',(1-a_i')T'_i+L_i'+R'_i+c_iS'_i)$ is klt. Since $a_i'>0$ and $\lim_{i\rightarrow+\infty}a_i'=0$, by \cite[Theorem 1.5]{HMX14}, possibly by passing to a subsequence, we may assume that $S'_i\not=0$.

Suppose that $\dim Z_i>0$ for infinitely many $i$. Possibly by passing to a subsequence, we may assume that $\dim Z_i>0$ for each $i$ and $\dim Z_i=\dim Z_j$ for each $i$ and $j$. Let $F_i$ be a general fiber of $X_i'\rightarrow Z_i$, and $B_{F_i}:=B'_i|_{F_i}$,  $L_{F_i}:=L'_i|_{F_i}$, $C_{F_i}:=C'_i|_{F_i}$, $R_{F_i}:=R'_i|_{F_i}$, $S_{F_i}:=S'_i|_{F_i}$, and $T_{F_i}:=T'_i|_{F_i}$. Then $(F_i,B_{F_i}=(1-a'_i)T_{F_i}+L_{F_i}+R_{F_i}+c_iS_{F_i})$ is klt and $K_{F_i}+B_{F_i}\equiv 0$. Since $T_i'$ dominates $Z_i$, $T_{F_i}\not=0$.  Since $a_i'>0$ and $\lim_{i\rightarrow+\infty}a_i'=0$, by \cite[Theorem 1.5]{HMX14}, possibly by passing to a subsequence, we may assume that $S_{F_i}\not=0$. Since $\dim F_i\leq\dim X_i'-1=d-2$, by \cite[Proposition 11.7]{HMX14}, $c\in N(d-3,\Ii)$. By Theorem \ref{thm: hmx accumulation point preliminaries}(3), $c\in\lct(d-2,\Ii)$, a contradiction. Thus possibly by passing to a subsequence, we may assume that $\dim Z_i=0$ for each $i$. In particular, $\rho(X_i')=1$ for each $i$.

\medskip

\noindent\textbf{Step 5}. We reduce our theorem to the case when $G_i'=0$ and $c_i':=\lct(X_i',T'_i+L_i'+R'_i;S'_i)\geq c_i$ in this step.

\smallskip

Since $(X_i',T'_i+L_i'+R'_i+cS'_i+G'_i)$ is lc, $(X_i',T'_i+L_i'+R'_i+cS'_i)$ is lc. We consider $c_i':=\lct(X_i',T'_i+L_i'+R'_i;S'_i)$, then $c_i'\geq c$. Possibly by passing to a subsequence, we may assume that either $c_i'\in [c,c_i)$ for each $i$ or  $c_i'\geq c_i$ for each $i$.

Suppose that $c_i'\in [c,c_i)$ for each $i$. Since $(X_i',(1-a'_i)T_i'+L_i'+R_i'+c_iS_i')$ is klt, all lc centers of $(X_i',T_i'+L_i'+R'_i+c_i'S'_i)$ are contained in $T_i'$, and there exists an lc center of $(X_i',T_i'+L_i'+R'_i+c_i'S'_i)$ which is contained in $T_i'\cap\Supp S_i'$. In particular, there exist general hyperplane sections $H_{i,1},\dots,H_{i,l_i}\subset X_i'$ for some integer $l_i\geq 0$ and $U_i:=X_i'\cap(\cap_{j=1}^{l_i}H_{i,j})$, such that
\begin{itemize}
    \item $(U_i,\Delta_i:=T_{U_i}+L_{U_i}+R_{U_i}+c_i'S_{U_i})$ is lc, where $T_{U_i}:=T_i'|_{U_i},L_{U_i}:=L_i'|_{U_i},R_{U_i}:=R'_i|_{U_i}$, and $S_{U_i}:=S'_i|_{U_i}$,
    \item all lc centers of $(U_i,\Delta_i)$ are contained in $T_{U_i}$,
    \item there exists an lc center of $(U_i,\Delta_i)$ which is contained in $T_{U_i}\cap\Supp S_{U_i}$, and
    \item all lc centers of $(U_i,\Delta_i)$ which are contained in $T_{U_i}\cap\Supp S_{U_i}$ have dimension $0$.
\end{itemize}
Possibly by passing to a subsequence, we may assume that $l_i=l_0$ is a constant. Let $g_i: W_i\rightarrow U_i$ be a dlt modification of $(U_i,\Delta_i)$, and $S_{W_i}$ the strict transform of $S_{U_i}$ on $W_i$. Note that $g_i^*S_{U_i}=S_{W_i}+F_i$ for some $F_i\geq 0$ such that $F_i\subset\Exc(g_i)$. We show that $\Supp S_{W_i}\cap F_i\not=\emptyset$. Suppose that $\Supp S_{W_i}\cap F_i=\emptyset$, then by the negativity lemma, $F_i=0$, hence $g_i$ is the identity morphism near $\Supp S_{U_i}$. Thus $(U_i,\Delta_i)$ is dlt near $\Supp S_{U_i}$. Since $c_i'<c_i<1$, $\lfloor \Delta_i\rfloor=T_{U_i}$, so $(U_i,\Delta_i)$ is plt near $\Supp S_{U_i}$. However, this is not possible since there exists an lc center of $(U_i,\Delta_i)$ which is contained in $T_{U_i}\cap\Supp S_{U_i}$.

Thus we can pick a $g_i$-exceptional prime divisor $E_i$ such that $E_i\cap\Supp S_{W_i}\not=\emptyset$ and $E_i\subset F_i$. Since $E_i$ is an lc place of $(U_i,\Delta_i)$, $\Center_{U_i}E_i$ is contained in $T_{U_i}$. Thus $V_i:=\Center_{U_i}E_i$ is contained in $T_{U_i}\cap\Supp S_{U_i}$, so $V_i$ is a point. 

We denote the sum of all $g_i$-exceptional prime divisors by $E_{g_i}$. We let $B_{W_i}:=(g_i^{-1})_*\Delta_i+E_{g_i}$ and
$K_{E_i}+B_{E_i}:=(K_{W_i}+B_{W_i})|_{E_i}$. Then $K_{E_i}+B_{E_i}\sim_{\mathbb R}0$ as $V_i$ is a point. Since $\dim E_i=d-2-l_0$ and $E_i\cap\Supp S_{W_i}\not=\emptyset$,
$$(E_i,B_{E_i})\in\mathfrak{N}(d-2-l_0,\Ii,c_i').$$
Since $V_i\in T_{U_i}$, $(E_i,B_{E_i})$ is not klt. By Lemma \ref{lem: coeff 1 imply low dimensional lct}, $c_i'\in N(d-3-l_0,\Ii)\subset N(d-3,\Ii)$. By Theorem \ref{thm: hmx accumulation point preliminaries}(1)(3)(4), $c\in N(d-3,\Ii)$, which is not possible.

Thus $c_i'\geq c_i$ for each $i$. Since $\rho(X_i')=1$, we may let $c_i''$ be the unique real number such that $K_{X_i'}+T_i'+L_i'+R_i'+c_i''S_i'\equiv 0.$ Since $K_{X_i'}+(1-a'_i)T_i'+L_i'+R_i'+c_iS_i'\equiv 0$ and $K_{X_i'}+T_i'+L_i'+R_i'+cS_i'+G_i'\equiv 0$, $c\leq c_i''<c_i\leq c_i'$. By Lemma \ref{lem: coeff 1 imply low dimensional lct}, $c_i''\in N(d-2,\Ii)$. Since $c\not\in N(d-3,\Ii)$, by Theorem \ref{thm: hmx accumulation point preliminaries}(4), possibly by passing to a subsequence, we may assume that $c_i''=c$ for each $i$. In particular, $G_i'=0$ for each $i$. 

\medskip

\noindent\textbf{Step 6}. In this step, we reduce our theorem to the case when $(X_i',T_i'+L_i'+R_i'+cS_i')$ is plt for each $i$.

\smallskip

Suppose that $(X_i',T_i'+L_i'+R_i'+cS_i')$ is not plt for infinitely many $i$. Possibly by passing to a subsequence, we may assume that $(X_i',T_i'+L_i'+R_i'+cS_i')$ is not plt for each $i$. Since $(X_i',(1-a_i')T_i'+L_i'+R_i'+cS_i')$ is klt, $(X_i',T_i'+L_i'+R_i'+cS_i')$ is not plt near $T_i'$. Since $\lct(X_i',T_i'+L_i'+R_i';S_i')=c_i'\geq c_i>c$, $(X_i',T_i'+L_i'+R_i'+cS_i')$ is plt near the generic point of each irreducible component of $T_i'\cap S_i'$. Moreover, since $(X_i',(1-a_i')T_i'+L'_i+R'_i+c_iS_i')$ is klt, any lc center of $(X_i',T_i'+L_i'+R_i'+cS_i')$ is contained in $T_i'$. Thus we may take a dlt modification $g_i: W_i\rightarrow X_i'$ of $(X_i',T_i'+L_i'+R_i'+cS_i')$ which is an isomorphism near the generic point of each irreducible component of $T_i'\cap S_i'$. We denote the sum of all $g_i$-exceptional prime divisors by $E_{g_i}$. We let $T_{W_i}$ and $S_{W_i}$ be the strict transforms of $T_i'$ and $S_i'$ on $W_i$ respectively. Since $T_i'\cap S_i'\not=\emptyset$ and  $g_i$ is an isomorphism near the generic point of each irreducible component of $T_i'\cap S_i'$, $T_{W_i}\cap S_{W_i}\not=\emptyset$. Let $B_{W_i}:=(g_i^{-1})_*(T_i'+L_i'+R'_i+cS'_i)+E_{g_i}$. Since $K_{X_i'}+T_i'+L_i'+R_i'+cS_i'\equiv 0$, $K_{W_i}+B_{W_i}\equiv 0$. Let $$K_{T_{W_i}}+B_{T_{W_i}}:=(K_{W_i}+B_{W_i})|_{T_{W_i}},$$
then $(T_{W_i},B_{T_{W_i}})\in\mathfrak{N}(d-2,\Ii,c)$ and $(T_{W_i},B_{T_{W_i}})$ is not klt since $(X_i',T_i'+L_i'+R_i'+cS_i')$ is not plt near $T_i'$. By Lemma \ref{lem: coeff 1 imply low dimensional lct}, $c\in N(d-3,\Ii)$, which is not possible.

Thus possibly by passing to a subsequence, we may assume that $(X_i',T_i'+L_i'+R_i'+cS_i')$ is plt for each $i$.  Since $G_i'=0$, $(X_i,L_i+R_i+cS_i+G_i)$ and $(X_i',T_i'+L_i'+R_i'+cS_i')$ are crepant. Since $(X_i,L_i+R_i+cS_i+G_i)$ is an $(N,\Ii_0)$-decomposable $\Rr$-complement of $(X_i,L_i+R_i+cS_i)$, $a(E_i,X_i,L_i+R_i+cS_i+G_i)\geq 2\epsilon_0$ for any prime divisor $E_i\not=T_i$ over $X_i$. 

\medskip

\noindent\textbf{Step 7}. We prove the case when $T_i$ is on $X_i$ for each $i$ in this step.

\smallskip

Suppose that $T_i$ is on $X_i$ for each $i$. Then $X_i=X_i'$ is $\epsilon_0$-klt, hence $X_i$ belongs to a bounded family by \cite[Theorem 1.1]{Bir21}. Moreover, since $1$ is the only possible accumulation point of $\Ii$, the coefficients of $L_i'$, $R_i'$, $S_i'$ belong to a finite set and $1-a_i'\in D(\Ii\cup\{c_i\})$. In particular, there exist a positive integer $M$ which does not depend on $i$, and very ample divisors $H_i$ on $X_i$, such that $-K_{X_i}\cdot H_i^{d-2}\leq M$. Since 
$$(K_{X_i}+T_i+L_i'+R_i'+cS_i')\cdot H_i^{d-2}=0$$
and $c>0$, $T_i\cdot H_i^{d-2}>0$ and $S_i'\cdot H_i^{d-2}$ belong to a finite set. Since $\rho(X_i)=1$, possibly by passing to a subsequence, we may assume that there exists a positive rational number $\lambda=\frac{p}{q}$, where $p$ and $q$ are coprime positive integers, such that $T_i\equiv \lambda S_i'$ for each $i$. Since
$$K_{X_i}+(1-a'_i)T_i+L_i'+R_i'+c_iS_i'\equiv 0\equiv K_{X_i}+T_i+L_i'+R_i'+cS_i',$$
we have 
$$a'_i\lambda S_i'\equiv a'_iT_i\equiv (c_i-c)S_i',$$
hence $c_i-c=a'_i\lambda$. Since $1-a_i'\in D(\Ii\cup\{c_i\})$, $$0<a'_i=\frac{1-\gamma_i-k_ic_i}{m_i}$$ for some $\gamma_i\in\Ii\backslash\{1\}$, $k_i\in\mathbb N$, and $m_i\in\mathbb N^+$. Since $c_i>c>0$, possibly by passing to a subsequence, we may assume that $k_i=k$ is a constant. Since $\Ii\subset\left\{1-\frac{b_j}{n}\mid j,n\in\mathbb N^+, 1\leq j\leq m\right\}$, possibly by passing to a subsequence, we have $a'_i=\frac{\frac{b_j}{n_i}-kc_i}{m_i}$
for $n_i\in\mathbb N^+$ and some fixed $j$. Thus
$$c_i=c+\frac{\frac{\lambda b_j}{n_i}-\lambda kc}{m_i+\lambda k}.$$
If $k=0$, then $c_i=c+\frac{\lambda b_j}{n_im_i}\in \left\{c+\frac{\lambda b_j}{n}\mid n\in\mathbb N^+\right\}$, hence $\{c_i\}_{i=1}^{+\infty}$ is standardized near $c$ and we are done. If $k>0$, then since $c_i>c$, possibly by passing to a subsequence, we may assume that $n_i=n_0$ is a constant, hence $c_i\in\left\{c+\frac{q(\frac{\lambda b_j}{n_0}-\lambda kc)}{n}\mid n\in\mathbb N^+\right\}$, so $\{c_i\}_{i=1}^{+\infty}$ is standardized near $c$ and we are done.

\medskip

\noindent\textbf{Step 8}. We conclude the proof in this step. By \textbf{Step 7}, possibly by passing to a subsequence, we may assume that $T_i$ is exceptional over $X_i$ for each $i$. Then $f_i: Y_i\rightarrow X_i$ is the birational contraction which only extracts $T_i$, and $f_i$ is an $\epsilon_0$-plt blow-up of $(X_i\ni x_i,L_i+R_i+cS_i)$, where $x_i$ is the generic point of $\Center_{X_i}T_i$. Moreover, the coefficients of $L_i,R_i$ belong to a finite set, and the coefficients of $S_i$ belong to a finite rational set. By Lemma \ref{lem: complicated plt blow up}, possibly by passing to a subsequence, there exist a positive real number $\alpha$ and a non-negative rational number $\beta$ such that $a_i'=\frac{\alpha-(c_i-c)\beta}{n_i}$ for some positive integer $n_i$.

Since $(X_i',T_i'+L_i'+R_i'+cS_i')$ is plt and $(X_i',T_i'+L_i'+R_i'+cS_i')$ is an $(N,\Ii_0)$-decomposable $\Rr$-complement of itself, $X_i'$ is $\epsilon_0$-klt. Thus $X_i'$ belongs to a bounded family by \cite[Theorem 1.1]{Bir21}, and there exist a positive integer $M$ and very ample divisors $H_i$ on $X_i'$ such that $-K_{X_i'}\cdot H_i^{d-2}\leq M$. Since 
$$(K_{X'_i}+T'_i+L_i'+R_i'+cS_i')\cdot H_i^{d-2}=0$$
and $c>0$, $T'_i\cdot H_i^{d-2}>0$ and $S_i'\cdot H_i^{d-2}$ belong to a finite set. Since $\rho(X_i')=1$, possibly by passing to a subsequence, we may assume that there exists a positive rational number $\lambda$ such that $T_i'\equiv\lambda S_i'$ for each $i$. Since $$K_{X_i'}+(1-a_i')T'_i+L_i'+R_i'+c_iS_i'\equiv 0\equiv K_{X'_i}+T'_i+L_i'+R_i'+cS_i',$$
we have 
$$a'_i\lambda S_i'\equiv a'_iT'_i\equiv (c_i-c)S_i',$$
hence $c_i-c=a'_i\lambda$. Thus
$$c_i-c=\frac{\alpha\lambda}{n_i+\beta\lambda}.$$
Let $\mu$ be a positive integer such that $\mu\beta\lambda\in\mathbb N$, then 
$$c_i=c+\frac{\mu\alpha\lambda}{\mu n_i+\mu\beta\lambda}\in\left\{c+\frac{\mu\alpha\lambda}{n}\mid n\in\mathbb N^+\right\}.$$
Thus $\{c_i\}_{i=1}^{+\infty}$ is standardized near $c$, and we are done.
\end{proof}

\begin{proof}[Proof of Theorem \ref{thm: standardized lcts}]
It follows from Theorem \ref{thm: standardization lcts with boundary}.
\end{proof}

\section{Standardization of threefold canonical thresholds}

\begin{proof}[Proof of Theorem \ref{thm: standardized threefold thresholds}]
We only need to show that $\ct(3)$ is standardized since $\ct(1)=\{0\}$ and $\ct(2)=\{\frac{1}{n}\mid n\in\mathbb N^+\}\cup\{0\}$ by \cite[Lemma 2.17]{HL22}. By \cite[Theorem 1.8]{HLL22} and \cite[Theorem 1.1]{Che22b}, we know $\partial\ct(3)=\{\frac{1}{k}\mid k\in\mathbb N^+,k\geq 2\}\cup\{0\}$. Thus we only need to show that $\ct(3)$ is standardized near $\frac{1}{k}$ for any integer $k\geq 2$. By Lemma \ref{lem: standardized iff subsequence standardized}, we only need to show that for any sequence $\{c_i\}_{i=1}^{+\infty}\subset\ct(3)$ such that $\lim_{i\rightarrow+\infty}c_i=\frac{1}{k}$, a subsequence of $\{c_i\}_{i=1}^{+\infty}$ is standardized near $\frac{1}{k}$. In the following, we fix $k$ and $\{c_i\}_{i=1}^{+\infty}\subset\ct(3)$. 

Possibly by passing to a subsequence, we may assume that $c_i\in (\frac{1}{k},\frac{1}{k-1})$ for each $i$, $c_i$ is strictly decreasing, and $c_i=\ct(X_i\ni x_i,0;B_i)$, where $X_i\ni x_i$ is a threefold terminal singularity and $B_i\geq 0$ is Weil divisor on $X_i$. Possibly by replacing $X_i$ with a $\mathbb Q$-factorialization, we may assume that $X_i$ is $\mathbb Q$-factorial for each $i$. By \cite[Theorem 4.8]{HLL22}, $0$ is the only accumulation point of canonical thresholds whose ambient variety is neither smooth, nor of $cA$-type or $cA/n$-type. Since $\lim_{i\rightarrow+\infty}c_i=\frac{1}{k}>0$, possibly by passing to a subsequence, we may assume that $X_i\ni x_i$ is either smooth, or a $cA$-type singularity, or a $cA/n_i$-type singularity for some positive integer $n_i$. By \cite[Propositions 2.1, 2.2]{Che22b}, we may assume that $X_i\ni x_i$ is a $cA/n_i$-type singularity for some positive integer $n_i$. By \cite[Lemma 5.10]{Che22a}, $n_i\leq 3k$, so possibly by passing to a subsequence, we may assume that $n=n_i$ is a constant.

By \cite[Claims 2.4, 2.5, 2.6]{Che22b} and \cite[Lemma 5.2]{Che22a}, we may assume that 
$$c_i=\ct(X_i,0;B_i)=\frac{a_i}{m_i},$$
and there exist positive integers $d_i$, non-negative integers $l_{2,i},l_{3,i},r_{1,i},r_{2,i}$, such that
\begin{itemize}
    \item $r_{1,i}+r_{2,i}=a_id_in$ and $r_{1,i}\leq r_{2,i}$ (\cite[Proof of Lemma 2.3, Line 8]{Che22b}),
\item if $a_i\nmid m_i$, then $m_i\geq\frac{r_{1,i}r_{2,i}}{d_in^2}$ (\cite[Lemma 5.2]{Che22a}),
    \item if $a_i\geq 6k^2$, then $d_in\leq 4k$ (\cite[Claim 2.4]{Che22b}),
    \item $\max\{l_{2,i},l_{3,i}\}<k$ and either $l_{2,i}>0$ or $l_{3,i}>0$ (\cite[Claim 2.5]{Che22b}), 
    \item $d_inl_{2,i}+l_{3,i}\geq k$ (\cite[Claim 2.6]{Che22b}), and
    and
    \item $$\frac{a_i}{m_i}\leq\frac{a_i-n}{(r_{2,i}-d_in^2)l_{2,i}+(a_i-n)l_{3,i}}$$
    (\cite[Claim 2.6]{Che22b}).
\end{itemize}
Possibly by passing to a subsequence, we may assume that $l_2:=l_{2,i}$ and $l_3:=l_{3,i}$ are constant integers. Since $\lim_{i\rightarrow+\infty}c_i=\frac{1}{k}$, $\lim_{i\rightarrow+\infty}a_i=\lim_{i\rightarrow+\infty}m_i=+\infty$. Therefore, possibly by passing to a subsequence, we have $d_in\leq 4k$. Possibly by passing to a subsequence, we may assume that $d:=d_i$ is a constant. Since $\frac{a_i}{m_i}\in (\frac{1}{k},\frac{1}{k-1})$, $a_i\nmid m_i$. Thus $m_i\geq\frac{r_{1,i}r_{2,i}}{dn^2}$, so
$$\frac{1}{k}<\frac{a_i}{m_i}=\frac{a_idn^2}{r_{1,i}r_{2,i}}=\frac{n}{r_{1,i}}+\frac{n}{r_{2,i}}\leq\frac{2n}{r_{1,i}},$$
so $r_{1,i}<2kn$. Possibly by passing to a sequence, we may assume that $r_1=r_{1,i}$ is a constant. Therefore,
$$\frac{a_i}{m_i}\leq\frac{a_i-n}{(r_{2}-dn^2)l_{2}+(a_i-n)l_{3}}=\frac{a_i-n}{(a_idn-r_1-dn^2)l_{2}+(a_i-n)l_{3}}=\frac{1-\frac{n}{a_i}}{dnl_2+l_3-\frac{(r_1+dn^2)l_2+nl_3}{a_i}}.$$
Since $\lim_{i\rightarrow+\infty}\frac{a_i}{m_i}=\frac{1}{k}$, $dnl_2+l_3=k$. Let $I:=(r_1+dn^2)l_2+nl_3$, then
$$\frac{1}{k}<\frac{a_i}{m_i}\leq\frac{1-\frac{n}{a_i}}{k-\frac{I}{a_i}}=\frac{a_i-n}{ka_i-I},$$
so
$$ka_i>m_i\geq ka_i-\frac{a_i}{a_i-n}(I-kn).$$
Since $\lim_{i\rightarrow+\infty}a_i=+\infty$, possibly by passing to a subsequence, we may assume that there exists a positive integer $I'$ such that $m_i=ka_i-I'$ for each $i$. Thus
$$c_i=\frac{a_i}{m_i}=\frac{a_i}{ka_i-I'}=\frac{1}{k}+\frac{I'}{k(ka_i-I')}\in\left\{\frac{1}{k}+\frac{I'}{m}\mid m\in\mathbb N^+\right\},$$
so $\{c_i\}_{i=1}^{+\infty}$ is standardized near $\frac{1}{k}$, and we are done.
\end{proof}


\begin{thebibliography}{99}


\bibitem[Ale93]{Ale93} V.~Alexeev, \textit{Two two--dimensional terminations}, Duke Math. J. \textbf{69} (1993), no. 3, 527--545.
 E. Miller (Ed.), Contemporary Mathematics \textbf{502} (2009), 1--4.

\bibitem[ADL23]{ADL23} K. Ascher, K. DeVleming, and Y. Liu, \textit{K-stability and birational models of moduli of quartic K3 surfaces}, Invent. Math. \textbf{232} (2023), 471--552.


	
\bibitem[Bir19]{Bir19} C. Birkar, \textit{Anti-pluricanonical systems on Fano varieties}, Ann. of Math. (2), \textbf{190} (2019), 345--463.
	


\bibitem[Bir21]{Bir21} C. Birkar, \textit{Singularities of linear systems and boundedness of Fano varieties}, Ann. of Math. (2) \textbf{193} (2021), 347--405.



\bibitem[BCHM10]{BCHM10}
C. Birkar, P. Cascini, C. D. Hacon and J. M\textsuperscript{c}Kernan, \textit{Existence of minimal models for varieties of log general type}, J. Amer. Math. Soc. \textbf{23} (2010), no. 2, 405--468.



\bibitem[CH21]{CH21} G. Chen and J. Han, \textit{Boundedness of $(\epsilon, n)$-complements for surfaces}, Adv. Math. \textbf{383} (2021), 107703, 40pp.

\bibitem[Che22a]{Che22a} J.-J. Chen, \textit{On threefold canonical thresholds}, Adv. Math. \textbf{404} (2022), Part B, 108447.

\bibitem[Che22b]{Che22b} J.-J. Chen, \textit{Accumulation points on 3-fold canonical thresholds}, arXiv:2202.06230.


\bibitem[Cor95]{Cor95} A. Corti, \textit{Factoring birational maps of 3-folds after Sarkisov}, J. Algebraic Geom. \textbf{4} (1995), 223--254.


\bibitem[HMX14]{HMX14} C. D. Hacon, J. M\textsuperscript{c}Kernan, and C. Xu, \textit{ACC for log canonical thresholds}, Ann. of Math. \textbf{180} (2014), no. 2, 523--571.



\bibitem[HLS19]{HLS19} J. Han, J. Liu, and V. V. Shokurov, \textit{ACC for minimal log discrepancies of exceptional singularities}, arXiv:1903.04338.

\bibitem[HLL22]{HLL22} J. Han, J. Liu, and Y. Luo, \textit{ACC for minimal log discrepancies of terminal threefolds}, arXiv:2202.05287.



\bibitem[HL22]{HL22} J. Han and Y. Luo, \textit{On boundedness of divisors computing minimal log discrepancies for surfaces}, J. Inst. Math. Jussieu. (2022), 1--24.

\bibitem[Jia21]{Jia21} C.~Jiang, \textit{A gap theorem for minimal log discrepancies of non-canonical singularities in dimension three}, J. Algebraic Geom. \textbf{30} (2021), 759--800.


\bibitem[Kaw88]{Kaw88} Y. Kawamata, \textit{Crepant blowing-up of $3$-dimensional canonical singularities and its application to degenerations of surfaces}, Ann. of Math. (2), \textbf{127} (1988), no. 1, 93--163.



\bibitem[Kol08]{Kol08} J. Koll\'ar, \textit{Which powers of holomorphic functions are integrable?}, arXiv:0805.0756.


\bibitem[KM98]{KM98} J. Koll\'{a}r and S. Mori, \textit{Birational geometry of algebraic varieties}, Cambridge Tracts in Math. \textbf{134} (1998), Cambridge Univ. Press.

\bibitem[Kud01]{Kud01} S. A. Kudryavtsev, \textit{Pure log terminal blow-ups}, Math. Notes, \textbf{69} (2001), no. 5, 814--819.



\bibitem[LX21]{LX21} J.~Liu and L.~Xiao, \textit{An optimal gap of minimal log discrepancies of threefold non-canonical singularities}, J. Pure Appl. Algebra \textbf{225} (2021), no. 9, 106674, 23 pp.

\bibitem[LL22]{LL22} J.~Liu and Y. Luo, \textit{Second largest accumulation point of minimal log discrepancies of threefolds}, arXiv:2207.04610.


\bibitem[MN18]{MN18} M.~Musta\c{t}\u{a} and Y.~Nakamura, \textit{A boundedness conjecture for minimal log discrepancies on a fixed germ}, Local and global methods in algebraic geometry, Contemp. Math. \textbf{712} (2018), 287--306.

\bibitem[Nak16]{Nak16} Y.~Nakamura, \textit{On minimal log discrepancies on varieties with fixed Gorenstein index}, Michigan Math. J. \textbf{65} (2016), no. 1, 165--187.

\bibitem[Pro00]{Pro00} Y.G. Prokhorov, \textit{Blow-ups of canonical singularities}, Algebra (Moscow, 1998) (2000): 301--318.



\bibitem[Pro18]{Pro18} Y.G.~Prokhorov, \textit{The rationality problem for conic bundles} (Russian), Uspekhi Mat. Nauk \textbf{73} (2018), no. 3 (441), 3--88; translation in Russian Math. Surveys \textbf{73} (2018), no. 3, 375--456. 

\bibitem[PS09]{PS09} Y.G. Prokhorov and V.V. Shokurov, \textit{Towards the second main theorem on complements}, J. Algebraic Geom., \textbf{18} (2009), no. 1, 151--199.


\bibitem[Rei87]{Rei87} M.~Reid, \textit{Young person’s guide to canonical singularities}, Algebraic geometry, Bowdoin 1985, Proc. Symp. Pure Math. \textbf{46} (1987), Part 1, 345--414.

\bibitem[Sho88]{Sho88} V.V.~Shokurov, {\it Problems about {F}ano varieties}, {Birational Geometry of Algebraic Varieties, Open Problems. The XXIIIrd International Symposium, Division of Mathematics, The Taniguchi Foundation}, 30--32, August 22--August 27, 1988.



\bibitem[Sho96]{Sho96} V.V. Shokurov, \textit{3-fold log models}, J. Math. Sci. \textbf{81} (1996), no. 3, 2667--2699.



\bibitem[Xu14]{Xu14} C. Xu, \textit{Finiteness of algebraic fundamental groups}, Compos. Math. \textbf{150} (2014), no. 3, 409--414.
\end{thebibliography}
\end{document}